\documentclass[12pt]{article}
\usepackage{amsmath,amsfonts,amssymb,amsthm,epsfig}%    epsfig added here
\newcommand{\mc}[1]{}
\newcommand{\bb}{\mathbb}
\newcommand{\gothic}{\mathfrak}

\newcommand{\integers}{{\bb Z}}

\newcommand{\reals}{{\bb R}}

\newcommand{\bbT}{{\bb T}}

\newlength{\figboxwidth}
\newcommand{\makefig}[3]{
        \begin{figure}[htb]
        \refstepcounter{figure}
        \label{#2}
        \begin{center}
                #3~\\
                \smallskip
                Figure \thefigure.  #1
        \end{center}
        \medskip
        \end{figure}
}
\newcommand{\makefignocenter}[3]{
        \begin{figure}[htb]
        \refstepcounter{figure}
        \label{#2}
        \begin{center}
                #3~\\
        \end{center}
                \smallskip

                Figure \thefigure.  #1
        \medskip
        \end{figure}
}
\newcommand{\appendixmode}{
        \setcounter{section}{0}
        \renewcommand{\thesection}{\Alph{section}}
}

\renewcommand{\bold}[1]{\medskip \noindent {\bf #1 }\nopagebreak}

\newcommand{\isom}{\cong}

\newcommand{\cross}{\times}
\newcommand{\st}{\;\: : \;\:}         %Such that

\newcommand{\zed}{\integers}

\newcommand{\area}{\operatorname{area}}

\def\@ifundefined#1#2#3%
  {\expandafter\ifx\csname#1\endcsname\relax#2\else#3\fi}

\@ifundefined{theoremstyle}{
}{
\theoremstyle{plain} %default
}
\newtheorem{theorem}{Theorem}[section]

\newtheorem{proposition}[theorem]{Proposition}
\newtheorem{lemma}[theorem]{Lemma}

\newtheorem{corollary}[theorem]{Corollary}

\@ifundefined{theoremstyle}{
}{
\theoremstyle{definition} %default
}

\newcommand{\cF}{{\cal F}}

\newcommand{\cH}{{\cal H}}

\newcommand{\cM}{{\cal M}}

\newcommand{\cP}{{\cal P}}

\newcommand{\cT}{{\cal T}}

\mathchardef\GG="321D
\newcommand{\gt}{{\gothic t}}

\newcommand{\semidirect}{\ltimes}

\title{Billiards in rectangles with barriers} 
\author{Alex Eskin\thanks{Research partially supported by 
                 NSF grant DMS-9704845, the  Sloan 
                Foundation and the Packard Foundation}, 
\ Howard Masur\thanks{Research partially supported by NSF
grant 9803497}
\ and Martin Schmoll 
\\ \\ }

\begin{document}
\maketitle

\section{Introduction}

Let $P=P(p/q, \alpha)$ denote the polygon whose boundary is the
union of the boundary of the unit square $\partial([0,1] \cross
[0,1])$ and the interval $\{p/q\} \cross [0, \alpha]$ (see
Figure~\ref{fig:1}). 
We consider ``billiard paths''; straight lines in the polygon,
that when hitting a side reflect with equal angles and that do
not pass through vertices. 
If a trajectory closes up, then it determines a parallel family of
closed trajectories of the same length. The boundary of such a
parallel
family consists of a pair of trajectories  
from the point $(p/q,\alpha)$
to itself. Let $v_1, v_2$ denote unit  vectors in the directions 
where the trajectories leave $(p/q,\alpha)$ and let $w_1, w_2$
denote unit vectors in the directions where the trajectories return to 
$(p/q,\alpha)$. Because of the fact that the polygon has right angles, 
the set $\{\pm v_1, \pm v_2, \pm w_1, \pm w_2 \}$ is stable under the
action of the 
Klein group (i.e the group generated by reflections in the horizontal
and vertical axes). 
For  each $T$ let $N_c(p/q,\alpha,T)$ denote the number of
parallel
families of nonprimitive  closed trajectories of Eulcidean length at most $T$.
(This allows trajectories to repeat themselves). 
We also consider ``generalized diagonals'', i.e. 
billiard trajectories  that join the point 
$(p/q,\alpha)$ to itself which are not on the boundary of a parallel
family; 
These can occur singly or in pairs, 
in the sense that there are two that are of the same length and the
set $\{\pm v_1, \pm w_1, \pm v_2, \pm w_2 \}$ of incoming and outgoing 
directions of the two trajectories is stable under the Klein group. 
Let $N_s^{(1)}(p/q,\alpha,T)$ denote the number of
such generalized diagonals of multiplicity $1$ of length at most
$T$ 
and $N_s^{(2)}(p/q,\alpha)$
the  number of multiplicity $2$.

\begin{theorem}
\label{theorem:asymp}
 There exist
constants $s_1(q),s_2(q)$ and $c(q)$ depending only on $q$ and
independent of $p$ and $\alpha$, 
such that whenever $\alpha$ is irrational, 
\begin{displaymath}
N_s^{(i)}(p/q, \alpha, T) \sim \frac{1}{4} s_i(q) \pi T^2, i=1,2
\end{displaymath} and 
\begin{displaymath}
N_c(p/q,\alpha,T)\sim \frac{1}{4} c(q)\pi T^2
\end{displaymath}
as $T \to \infty$. 
\end{theorem}

In the classical integrable case of billiards in the unit square,
the asymptotics are well-known, since a parallel family of length
at most $T$ corresponds to an integer lattice point $(2m,2n)$
with $4(m^2+n^2)\leq T^2$.  The asymptotics are $\frac{1}{4}\pi
T^2$.
The equilateral triangle, the isoceceles right triangle and the
$\pi/3$ right triangle are also integrable and thus have easily
identified asymptotics. Veech \cite {Vee} found  
a surprising class of nonintegrable of billiards 
 for which one can find precise asymptotics for the number of
cylinders
and saddle connections. 
  One such set of examples are billiards in right triangles with
 angle $\pi/n, n\geq 5$.
These billiards have the property that a certain affine
automorphism group of the billiard is a lattice in 
$SL(2,\reals)$.   Veech showed that whenever that is the case,
 the counting problems have exact quadratic asymptotics. Another
proof is given in  \cite {G-J2}.
  The billiards 
$P(p/q,\alpha)$ are Veech precisely when $\alpha$ is rational. 
Thus our theorem provides, for $\alpha$ irrational, a new class of
billiards for which precise asymptotics are known.  

\begin{theorem}
\label{theorem:constants}
We have $c(2) = 9/2$, $s_1(2) = 0$, $s_2(2) = 2$. For $q \ge 3$, 
\begin{displaymath}
c(q) = \frac{10 q - 11}{2 q - 2}
\end{displaymath}
\begin{displaymath}
s_1(q) = \frac{27(q - 2)}{8(q - 1)}
\end{displaymath}
\begin{displaymath}
s_2(q) = \frac{ 5 q + 6}{ 8(q-1)}
\end{displaymath}
\end{theorem}

\mc{BEGIN CHANGE}
We note that in the case where $q = 2$, both
Theorem~\ref{theorem:asymp} and Theorem~\ref{theorem:constants} 
were previously observed in \cite{Gutkin:Judge:private} using 
a more elementary argument. 
\mc{END CHANGE}

In this paper, we prove alternative formulas for the constants
$c(q)$, $s_1(q)$ and $s_2(q)$, see Proposition~\ref{prop:value:c:d},
Proposition~\ref{prop:value:s1} and
Proposition~\ref{prop:value:s2:d}. The proofs of the identities
which show that the alternative expressions are in fact
equal to the constants in Theorem~\ref{theorem:constants} will appear
in \cite{EOZ}. 
\medskip

%Explicit expressions for $N_q^P(2)$, $N_q^P(1,1)$, $s_2(q)$ and
%$c(q)$ are 
%given in Lemma~\ref{lemma:NP:2}, Lemma~\ref{lemma:NPd:11} and
%Proposition~\ref{prop:value:c:d}. 
%\medskip

Let $\alpha$ be a partition of $2 g
- 2$ (i.e. a representation of $2g - 2$ as a sum of positive
integers). 
Let $\cH(\alpha)$ denote the space of pairs $(M,\omega)$
where $M$ is a Riemann surface of genus $g$, and $\omega$ is a
holomorphic differential such that the orders of its zeroes is
given
by $\alpha$. This space is called a stratum; this term is
justified by
the fact that the space of all Abelian differentials on Riemann
surfaces of genus $g$ is stratified by the spaces $\cH(\alpha)$, 
as $\alpha$ varies over the partitions of $2g-2$. In particular, 
the stratum corresponding to the partition $(1,\dots,1)$ is
called the 
principal stratum: it corresponds to holomorphic differentials
with simple
zeroes.

Given a pair $(M,\omega)$ and a point $p \in M$ such that
$\omega(p)
\ne 0$, there exists a local coordinate $z$ near $p$ such that
$\omega 
= dz$. Such a  local coordinate is unique up to $z \to
z+c$. Thus $|dz|^2$ is a flat metric on $M$;  this metric
develops
conical singularities at the zeroes of $\omega$. An easy
calculation
shows that 
the total   angle at a zero of order $k$ is $2 \pi (k+1)$. 

Furthermore using 
such a coordinate near every nonsingular point
determines a structure of a ``translation surface'' on $M$,
namely an
atlas of coordinate charts which cover the surface away from the 
singularities, such that the transition functions are
translations $z
\to z + c$. Hence we may visualize $S = (M,\omega)$ as a union of
polygons contained in $\reals^2$ glued along parallel sides, 
such that each side is glued to exactly one other, and the angle
excess at each vertex is an integer multiple of $2\pi$. Geodesics
on $S$ are straight lines away from the singularities. We also 
define a {\em saddle connection} to be a straight line segment
which
begins and ends at a singularity. 

In genus $2$, there are 2 strata, namely the principal stratum
$\cH(1,1)$, consisting of holomorphic differentials with 2 simple
zeroes, and $\cH(2)$, consisting of holomorphic differentials
with one 
double zero. 
We take four copies of $P = P(p/q, \alpha)$ which are images of $P$
under reflection in the two coordinate axes and reflection in the origin
and glue them
as in Figure~\ref{fig:1}. Then we make the following identifications
on the union of the $4$ copies: we identify the top to the bottom, the 
left side to the right side. The interval $p/q \cross [0,\alpha]$
gives rise to two vertical double lines. We glue the left side of the
right line to the right side of the left line, and the right side of
the right line to the left side of the left line. 
We get a surface $S = 
S(p/q, \alpha)$ of area $4$, which belongs to $\cH(1,1)$, with the
$2$ zeroes located at the endpoints of the vertical lines. 
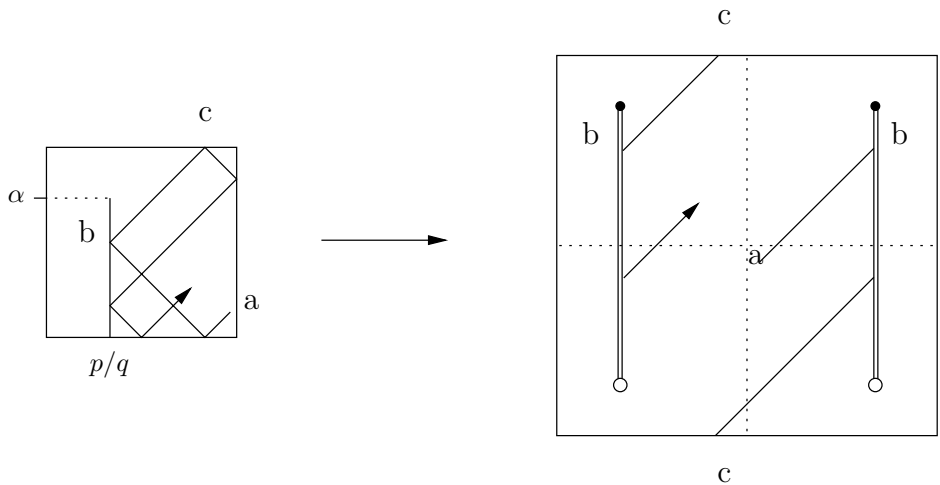
\makefig{On the invariant surfaces $S(p/q,\alpha)$ of the wall billiard \\ 
$P(p/q,\alpha)$  the trajectories are straight lines.}{fig:1}{ 
\begin{picture}(0,0)%
\includegraphics{developwall.pstex}%
\end{picture}%
\setlength{\unitlength}{4144sp}%
\begingroup\makeatletter\ifx\SetFigFont\undefined%
\gdef\SetFigFont#1#2#3#4#5{%
  \reset@font\fontsize{#1}{#2pt}%
  \fontfamily{#3}\fontseries{#4}\fontshape{#5}%
  \selectfont}%
\fi\endgroup%
\begin{picture}(5562,2835)(661,-2851)
\put(1151,-2183){\makebox(0,0)[lb]{\smash{\SetFigFont{10}{12.0}{\rmdefault}{\mddefault}{\updefault}
$p/q$
}}}
\put(661,-1171){\makebox(0,0)[lb]{\smash{\SetFigFont{10}{12.0}{\rmdefault}{\mddefault}{\itdefault}
$\alpha$
}}}
\put(2071,-1816){\makebox(0,0)[lb]{\smash{\SetFigFont{12}{14.4}{\rmdefault}{\mddefault}{\updefault}
a
}}}
\put(5086,-1546){\makebox(0,0)[lb]{\smash{\SetFigFont{12}{14.4}{\rmdefault}{\mddefault}{\updefault}
a
}}}
\put(5941,-826){\makebox(0,0)[lb]{\smash{\SetFigFont{12}{14.4}{\rmdefault}{\mddefault}{\updefault}
b
}}}
\put(4096,-826){\makebox(0,0)[lb]{\smash{\SetFigFont{12}{14.4}{\rmdefault}{\mddefault}{\updefault}
b
}}}
\put(1081,-1411){\makebox(0,0)[lb]{\smash{\SetFigFont{12}{14.4}{\rmdefault}{\mddefault}{\updefault}
b
}}}
\put(1801,-691){\makebox(0,0)[lb]{\smash{\SetFigFont{12}{14.4}{\rmdefault}{\mddefault}{\updefault}
c
}}}
\put(4906,-106){\makebox(0,0)[lb]{\smash{\SetFigFont{12}{14.4}{\rmdefault}{\mddefault}{\updefault}
c
}}}
\put(4906,-2851){\makebox(0,0)[lb]{\smash{\SetFigFont{12}{14.4}{\rmdefault}{\mddefault}{\updefault}
c
}}}
\end{picture}
}

A billiard trajectory $\lambda$ on $P(p/q,\alpha)$ gives rise to a
straight line  in a direction $\theta$ on $S  = S(p/q,\alpha)$; if
the line is in the upper left hand square of $S$ then the trajectory
$\lambda$ is moving in the direction $\theta$ on $P = P(p/q,\alpha)$;
if it is in the upper right hand square it is moving in the direction 
$\pi - \theta$; if it is in the lower right square it is in the
direction $\pi +\theta$ and if it is in the lower left square it is $2 \pi 
- \theta$.  (see Figure~\ref{fig:1}).

 A saddle
 connection on $S$ that returns to the zero 
from which it left, returns at  an angle of $\pi$ and hence is on
the boundary of a parallel family of closed geodesics. Thus a 
generalized diagonal on $P$ 
that is not the boundary of a parallel family corresponds to a
saddle connection joining distinct zeroes. In particular,
$N_s^{(i)}(p/q,
\alpha, T) = 
N_s^{(i)}(S(p/q, \alpha), T)$, and $N_c(p/q,\alpha, T) =
N_c(S(p/q,\alpha,T))$, where $N_s^{(i)}(S, T)$ (resp. $N_c(S,T)$)
denotes the 
number of saddle connections of multiplicity $i=1,2$
 (resp. cylinders of closed geodesics
which do not pass through singular points) on $S$ of length at
most $T$. (By a saddle connection of multiplicity $2$ 
we mean a pair of saddle connections of the same length and direction, 
connecting the same zeroes). 

\relax From \cite{Eskin:Masur:ae}, we have that there exist
constants $s_1(1,1),s_2(1,1)$ and $c(1,1)$
such that  
for almost all surfaces $S \in \cH(1,1)$, we have 
\begin{displaymath}
N_s^{(i)}(S,T) \sim s_i(1,1) \pi \frac{T^2}{\area(S)}
\end{displaymath}  and 
\begin{displaymath}
N_c(S,T)\sim c(1,1)\pi \frac{T^2}{\area(S)}
\end{displaymath}
as $T\to\infty$. These constants have been evaluated in 
\cite{Eskin:Masur:Zorich:SV}; we have $s_1(1,1) = 27/8,
s_2(1,1)=5/8$ 
and $c(1,1) = 5$. We emphasize that the constants $s_i(q)$ and $c(q)$
of Theorem~\ref{theorem:asymp} differ from the ``generic'' constants
$s_i(1,1)$ and 
$c(1,1)$. However, we have the following: 
\begin{theorem}
\label{theorem:constants:converge}
\begin{displaymath}
\lim_{q \to \infty}
s_1(q)=\frac{27}{8} =s_1(1,1)
\end{displaymath}
\begin{displaymath}
\lim_{q\to\infty}
s_2(q)=\frac{5}{8} = s_2(1,1)
\end{displaymath}
\begin{displaymath}
\lim_{q\to\infty}c(q)=5 = c(1,1)
\end{displaymath}
\end{theorem} 
The above is consistent with the behavior for unipotent flows 
on homogeneous spaces \cite{Dani:Margulis:distribution}. 

\bold{Remark:} Theorem~\ref{theorem:asymp} is a straightforward
application of Ratner's theorem and associated techniques. The
results 
about the evaluation of the constant
(i.e. Theorem~\ref{theorem:constants} or
Proposition~\ref{prop:value:c:d}, Proposition~\ref{prop:value:s1} and
Proposition~\ref{prop:value:s2:d}) require also certain results
about connectedness of the space of primitive torus covers
(Theorem~\ref{theorem:connected} below). Even though
Theorem~\ref{theorem:constants:converge} follows immediately from 
Theorem~\ref{theorem:constants}, in order to make the paper self contained
we derive it in
Appendix~A from the alternative expressions
Propositions~\ref{prop:value:s1}-\ref{prop:value:s2:d}.

%We note that $c(1,1)$ and $s_i(1,1)$ are themselves evaluated in
%terms
%of a ratio of volumes, and the volumes are also computed by
%counting
%torus covers, see \cite{Eskin:Masur:Zorich:SV} and
%\cite{Eskin:Okounkov:volumes}.

\section{Proof of Theorem~\ref{theorem:asymp}}
The essential feature of the surface $S(p/q, \alpha)$ which we
use is
that $S = (M, \omega)$ is a {\em $q$-fold torus cover} i.e. there
exists a holomorphic map $\pi: M \to \gt$, where $\gt$ is a
torus, and
such that $\omega$ is the pullback of the standard differential
$dz$
on $\gt$.  The map $\pi$ is a $q$-fold covering map, with
branching at
the zeroes of $\omega$.
Namely, we can divide the surface $S(p/q,\alpha)$ into $q$ vertical
strips of width $2/q$ and height $2$ and map each to a torus of the
same width and height.

In this section, we work in a
more general setup. Suppose $\cH(\beta)$ is a stratum (of Abelian
differentials), and suppose $S \in \cH(\beta)$ is a torus cover.
We
prove that $N_s^{(i)}(S,T) \sim s_i T^2$ and $N_c(S,T)\sim cT^2$
for
certain
constants $s_i,c$.

Since a torus covers itself, the torus covered by $S=(M,\omega)$
is
not unique.  We say that a torus cover $\pi': M \to\gt'$ factors
through $\pi: M \to \gt$ if there is an unbranched cover $\pi'':
\gt
\to \gt'$ such that $\pi' = \pi'' \circ \pi$.  Let $z_1, \dots,
z_m$
denote the zeroes of $\omega$, with $m = \ell(\beta)$.  We say
that
the covering $\pi: M \to \gt$ is {\em unstable} if there exists a
covering $\pi': M \to \gt'$ which factors through $\pi$ such that
the
number of distinct points in the set $\{\pi'(z_1), \dots,
\pi'(z_m) \}
\subset \gt'$ is strictly smaller then the number of distinct
points
in the set $\{ \pi(z_1), \dots, \pi(z_m) \} \subset \gt$. 
It is easy
to see that $S(p/q,\alpha)$ is a stable $q$-fold torus cover if
and
only if $\alpha$ is irrational. 

\mc{BEGIN CHANGE}
Let $\cM_q(\beta) \subset \cH(\beta)$ denote the surfaces
in $\cH(\beta)$ which have area $1$ and are $q$-fold torus covers. 
\mc{END CHANGE}

\begin{lemma}
\label{lemma:closed:subset}
The set $\cM_q(\beta)$ is a closed subset of $\cH(\beta)$
which is invariant under the action of $SL(2,\reals)$ on
$\cH(\beta)$. 
\end{lemma}

\bold{Proof.} 
Suppose $\omega_n$ is a sequence in $\cM_q(\beta)$ such there is
a sequence of tori $\gt_n$ so that $\omega_n$ is the pull-back of
$dz$ under a covering map 
$\pi:M\to\gt_n$. Suppose $\omega_n$  converges to $\omega_0$.  
The sequence $\gt_n$ cannot converge to the cusp; otherwise as a
$q$ fold cover$\omega_n$ would contain a curve whose length
approaches $0$. Thus passing to a subsequence we can assume
$\gt_n$ converges to $\gt_0$. The  holomorphic maps
$\pi:M\to\gt_n$ must converge to a holomorphic map which is a $q$
fold covering of $\omega_0$ over $\gt_0$. Hence $\cM_q(\beta)$ is
closed. Since $SL(2,\reals)$
acts linearly on tori as well as on flat structures, it is clear
that the $SL(2,\reals)$ orbit of a torus cover again consists of
torus covers.  
\qed\medskip

For $k > 0$ let $\cT^k$ denote the moduli space of tori with $k$
marked points; we always assume that one of the marked points is
at
the origin. Then, 
\begin{displaymath}
\cT^k \isom SL(2,\reals) \semidirect (\reals^2)^{k-1} / SL(2,
\zed)
\semidirect (\zed^2)^{k-1} \equiv G / \Gamma
\end{displaymath}
\begin{lemma}
\label{lemma:structure:Mq}
There is a covering map $p: \cM_q(\beta) \to \cT^k$, 
where $k$ denotes the number of parts of $\beta$ (i.e. the number
of
zeroes). The covering map commutes with the $SL(2,\reals)$
action.
(The action on $\cT^k = G/\Gamma$ is by left multiplication which
makes sense since $SL(2,\reals) \subset G$). 
\end{lemma}
\bold{Proof.} The map $p$ just sends $S = (M, \omega)\in
\cM_q^{(\beta)}$ to $(\pi(M), \pi(z_1) , \dots, \pi(z_k))$, where
$z_1, \dots, z_k$ are the zeroes of $\omega$.
\qed\medskip

For
any continuous compactly supported function $f: \reals^2 \to
\reals$, 
let $\hat{f}: \cH(\beta) \to \reals$
denote the Siegel-Veech transform
\begin{displaymath}
\hat{f}(S) = \sum_{v \in V(S)} f(v)
\end{displaymath}
where $V(S) \subset \reals^2$ denotes either the vectors
associated with
saddle connections of multiplicity $i$, or  
(cylinders of) periodic geodesics on $S$ which do not pass
through
singularities (see \cite{Eskin:Masur:ae} for details on the
notation). If
$V(\cdot)$ denotes the vectors associated with saddle
connections of multiplicity $i$,
then $N_s^{(i)}(S,T) = |V(S) \cap B(T)|$, where $B(T)$ is the
ball of
radius 
$T$ centered at the origin. If $V(\cdot)$ denotes the vectors
associated with cylinders of periodic geodesics, then $N_c(S,T) =
|V(S) \cap B(T)|$. 

Let $\mu$ denote the
Haar measure on $\cT^k$; then $\tilde\mu=p^{-1}(\mu)$ is a smooth
invariant measure for the $SL(2,\reals)$ action on
$\cM_q(\beta)$.

\begin{theorem}
\label{theorem:general:asymp}
Let $S \in \cM_q(\beta)$ denote a stable $q$-fold torus cover.
Then, 
there exists a constant $\kappa(S)$ such that as $T \to \infty$, 
\begin{displaymath}
|V(S) \cap B(T) | \sim \pi \kappa(S) T^2, 
\end{displaymath}
where $B(T)$ is the ball of radius $T$ (cf. \cite[Theorem
2.1]{Eskin:Masur:ae}). The constant $\kappa(S)$ depends only on
the
connected component $\cM(S)$ of $\cM_q(\beta)$ containing $S$. In
fact, $\kappa(S)$ is given by the following Siegel-Veech formula:
for any
continuous compactly supported $f: \reals^2 \to \reals$, 
\begin{displaymath}
\frac{1}{\tilde{\mu}(\cM(S))}\int_{\cM(S)} \hat{f} \,
d\tilde{\mu} =
\kappa(S) \int_{\reals^2} f
\end{displaymath}
\end{theorem}

\mc{BEGIN CHANGE}
We note that the surface $S(p/q,\alpha)$
which one obtains by gluing together four copies of the polygon
$P(p/q, \alpha)$ is a stable $q$-fold torus cover (and after rescaling 
by a factor of $1/2$ will have area $1$). 
Now Theorem~\ref{theorem:asymp} will follow from
Theorem~\ref{theorem:general:asymp}, the fact that these surfaces
are primitive covers (see \S\ref{sec:primitive}), and
Theorem~\ref{theorem:connected} which says that
the space of primitive covers is connected.
\medskip
\mc{END CHANGE}

We now begin the proof of Theorem~\ref{theorem:general:asymp}. 
As in \cite{Eskin:Masur:ae}, let $r_\theta = \begin{pmatrix} \cos
\theta & \sin
\theta \\ - \sin \theta & \cos \theta \end{pmatrix}$ and $a_t =
\begin{pmatrix} e^{t} & 0 \\ 0 & e^{-t} \end{pmatrix}$.
\relax From Proposition 3.2 of \cite{Eskin:Masur:ae},
it is enough to show that
\begin{equation}
\label{eq:main}
\lim_{t \to \infty} \frac{1}{2 \pi} \int_0^{2 \pi} \hat{f}( a_t
r_\theta S) \, d\theta = 
\frac{1}{\tilde{\mu}(\cM(S))}\int_{\cM(S)} \hat{f} \,
d\tilde{\mu}
\end{equation}
Now in view of the assumption that $S \in \cM_q(\beta)$ and
Lemma~\ref{lemma:closed:subset} the entire integral takes place
in
$\cM_q(\beta)$, in fact in the connected component $\cM(S)$ of
$\cM_q(\beta)$ containing $S$. Now in view of
Lemma~\ref{lemma:structure:Mq}, $\cM(S)$ is a finite cover of
$\cT^k$.
We again denote the covering map by $p$. Also we continue to
denote
the restriction of $\tilde{\mu}$ to $\cM(S)$ by $\tilde{\mu}$. 

The first step in the proof of
Theorem~\ref{theorem:general:asymp} is
the following:
\begin{lemma}
\label{lemma:measures:converge}
Let $a_t$, $r_\theta$, $S$ and $\cM(S)$ be as in
Theorem~\ref{theorem:general:asymp}. 
Let $K = \{ r_\theta \st 0 \le \theta \le 2 \pi \}$, so that $K
\isom 
SO(2)$. Let $\nu$ denote the $K$-invariant measure supported on
the
set $K S$. Then in the weak-* topology on $\cM(S)$, 
\begin{displaymath}
\lim_{t\to\infty}a_t \nu=\frac{1}{\tilde{\mu}(\cM(S))}\tilde\mu
\end{displaymath}
\end{lemma}
\bold{Proof.} 
Since the action of $SL(2,\reals)$ on $\cT^k$ is ergodic, and
$\cM(S)$ is connected, the action on $\cM(S)$ is ergodic.  

\mc{BEGIN CHANGE}
By \cite[Corollary~5.3]{Eskin:Masur:ae}, \mc{END CHANGE} 
for any sequence $t_i$, there exists a subsequence, again denoted
$t_i$ such that the measures $a_{t_i}\nu$ converge weakly to a
probability measure.  Hence it is enough to show that whenever
$a_{t_i}\nu$ converges to a measure $\nu_\infty$, we have
$\nu_\infty=\lambda \tilde \mu$, where $\lambda =
1/\tilde{\mu}(\cM(S))$.  

Consider the orbit of $SL(2,\reals)$ through $p(S)\in \cT^k\equiv
G/\Gamma$.  By Ratner's theorem  the orbit closure is itself the
orbit of a connected Lie group $F$ such that $SL(2,\reals)\subset
F\subset G$.  It is easy to check  that the only connected
subgroups $F$ of $G$ containing $SL(2,\reals)$ are of the form
$SL(2,\reals)\semidirect (\reals^2)^j$ for $1\leq j\leq k-1$.  Now the
stability assumption forces $F=G$, i.e. that the $SL(2,\reals)$
orbit through $p(S)$ is dense in $G/\Gamma$.

Since $p$ commutes with the $SL(2,\reals)$ action, it follows
from the
assumption that $a_{t_i} \nu \to \nu_\infty$ that $a_{t_i} p(\nu)
\to
p(\nu_\infty)$; the latter convergence takes place on $\cT^k$. 
Now by
\cite[Corollary 1.2]{Shah:SL2} (which uses Ratner's theorem), we
have
$p(\nu_\infty) = c \mu$, where $c \in \reals$ and $\mu$ is the
Haar
measure on $G/\Gamma$. Then $\nu_\infty$ is a measure on $\cM(S)$
which is absolutely continuous with respect to $\tilde \mu =
p^{-1}(\mu)$.  Also $\nu_\infty$ is $SL(2,\reals)$ invariant. 
Then by
the Radon-Nykodim theorem and the ergodicity of the action,
$\nu_\infty = \lambda \tilde{\mu}$, where $\lambda \in \reals$. 
Since
$\nu_\infty$ is a probability measure, $\lambda =
1/\tilde{\mu}(\cM(S))$.  
\qed\medskip

\bold{Proof of Theorem~\ref{theorem:general:asymp}.}  It is
enough to
establish (\ref{eq:main}).  However, the equation (\ref{eq:main})
does
not immediately follow from Lemma~\ref{lemma:measures:converge}
because the function $\hat f$ is not bounded (it is continuous on an
open dense set of full measure). 

As in \cite{Eskin:Masur:ae}, 
let $\ell(S)$ denote the length of the shortest saddle
connection on $S$.  For $\epsilon>0$, let $h_\epsilon$ be a
smooth
function such that $h_\epsilon(S)=1$ if $\ell(S)>\epsilon$ and
$h_\epsilon(S)=0$ if $\ell(S)<\epsilon/2$. Then for any
$\eta>\delta>0$, 
\begin{equation}
\label{eq:hat:f:one:minus:h}
\hat f(S)(1-h_\epsilon(S))\leq
\frac{C(\delta)}{\ell(S)^{1+\delta}}(1-h_\epsilon(S))\leq
\epsilon^{\eta-\delta}\frac{C(\delta)}{\ell(S)^{1+\eta}},
\end{equation}
where the left inequality is by \cite[Theorem 5.1]{Eskin:Masur:ae}
and
the right is by the fact that $1-h_\epsilon$ is supported on the
set
where $\ell(\cdot)<\epsilon$.  Let $A_t$ denote the averaging
functional $A_t(\phi)=\frac{1}{2\pi}\int_0^{2\pi}\phi(a_tr_\theta
S)d\theta$.  Then, by (\ref{eq:hat:f:one:minus:h}),
\begin{displaymath}
A_t(\hat fh_\epsilon)\leq A_t(\hat f)=A_t(\hat
fh_\epsilon)+A_t(\hat f(1-h_\epsilon))\leq A_t(\hat f
h_\epsilon)+\epsilon^{\eta-\delta}C(\delta)A_t(\ell(\cdot)^{-1-
\eta}).
\end{displaymath}
By \cite[Theorem 5.2]{Eskin:Masur:ae} $A_t(\ell(\cdot)^{-1-\eta})$
is uniformly
bounded by a constant $C_1(\eta)$ 
as $t\to\infty$.  By Lemma~\ref{lemma:measures:converge} and
since $\hat fh_\epsilon$ is continuous  and compactly
supported,
$A_t(\hat fh_\epsilon)$ converges to
$\lambda\int_{\cM(S)}  \hat f h_\epsilon d\tilde \mu$, where
$\lambda =
\frac{1}{\tilde{\mu}(\cM(S))}$. 
Hence
\begin{displaymath}
\lambda \int_{\cM(S)}\hat fh_\epsilon
d\tilde\mu\leq \liminf_{t \to \infty} A_t(\hat 
f)\leq \limsup_{t \to \infty} A_t(\hat f)\leq
\lambda \int_{\cM(S)}\hat fh_\epsilon 
d\tilde\mu+C(\delta)C_1(\eta)\epsilon^{\eta-\delta}
\end{displaymath}
Letting $\epsilon\to 0$ we get 
\begin{displaymath}
\lim_{t\to\infty}\frac{1}{2\pi}\int_0^{2\pi}\hat f(a_tr_\theta
S)d\theta=\lambda \int_{\cM(S)} \hat f d\tilde\mu
\end{displaymath} 
as required.
\qed\medskip

\section{Cylinder decompositions}
In this section we introduce coordinates on the spaces
$\cM_d(1,1)$
and $\cM_d(2)$. These coordinates will be used throughout the
paper. 
 
\subsection{The space $\cM_d(1,1)$.}
\label{sec:Md11}
Recall that $S = (M,\omega) \in \cM_d(1,1)$ if and only if
$S$ has area $1$ and 
there exists a $d$-fold holomorphic covering map $\pi: M \to
\gt$,
where $\gt$ is a torus, and $\omega$ is the pullback of the
differential $dz$ on $\gt$.  The covering $\pi$ is branched of
order $2$ 
at the zeroes
$z_1$ and $z_2$ of $\omega$. Note that if $(M,\omega)$ has area
$1$, 
$(\gt, dz)$ has area $1/d$. 

Let $\cT^2$ denote the moduli space of tori with two marked
points. 
Then we have the covering map $p: \cM_d(1,1) \to \cT^2$, sending 
$(M,\omega)$ to $(\pi(M),\pi(z_1),\pi(z_2))$, (see
Lemma~\ref{lemma:structure:Mq}). 

We may write $\gt = \reals^2/\Delta$, where $\Delta \subset
\reals^2$
is a lattice. Let $(v_1, v_2)$ be a reduced basis for the
lattice;
i.e. $(v_1, v_2)$ is a basis for $\Delta$ such that $(v_1, v_2)
\in
\cF$, where the fundamental domain $\cF$ is given by
\begin{displaymath}
\cF = \{ (v_1, v_2 ) \in \reals^2 \cross \reals^2 \st \|v_1\| \le
\|v_2\|,\  \langle v_1,v_2 \rangle \le \tfrac{1}{2} \langle
v_1,v_1
\rangle ,\ |v_1 \cross v_2 | = 1/d \}. 
\end{displaymath}
We always assume that $\pi(z_1)$ is the origin in
$\reals^2/\Delta$. 
We may write $\pi(z_2) = \delta_1 v_1 + \delta_2 v_2$, where $0
\le
\delta_1 < 1$, and $0 \le \delta_2 < 1$. Hence the tuple $((v_1,
v_2),
\delta_1, \delta_2) \in \cF \cross [0,1) \cross [0,1)$ can be
used as
a coordinate system on $\cT^2$.

Let $(\cT')^2$ denote the subset of $\cT^2$ with $\delta_2 \ne
0$, 
and let $\cM'_d(1,1)$ denote the surfaces in $\cM_d(1,1)$ which
project
under $p$ to $(\cT')^2$. It is clear that $\cM'_d(1,1)$ is an
open
and dense subset of $\cM_d(1,1)$. Our coordinate system will only
be
valid on $\cM'_d(1,1)$. 

We define the direction of $v_1$ to be 
{\em horizontal} and the direction of $v_2$ to be {\em vertical}.
In using these terms we are not assuming that the directions are
perpendicular.  However, all figures are drawn in the case where the 
horizontal and vertical directions are perpendicular. 
Recall that the  Abelian differential $\omega$ is the pullback
$\pi^{
  \ast}d z$ of the canonical differential $d z$ 
on $\gt$. Consequently any closed direction on the base surface
is a closed direction on $M$. In particular, the horizontal (and
vertical) trajectories through any point of $M$ are closed.

We now assume $(M,w) \in \cM'_d(1,1)$.  By taking the preimages
(under $\pi$) of the two horizontal closed geodesics through
$\pi(z_1)$ and $\pi(z_2)$ one obtains four saddle connections on
$(M,
\omega)$, with two connecting each $z_i$ with itself. Because the
horizontal direction on the base torus is closed, any horizontal
trajectory away from the four saddle connections, is a closed
geodesic
on $M$. The closed geodesics occur in parallel families of the
same
length. Topologically they are cylinders whose boundary consists
of
saddle connections. Note that by construction, the length of
these
geodesics has to be an integer multiple of $\|v_1\|$.

\begin{lemma}
\label{lemma:decomp}
Suppose $(M,\omega) \in \cM'_d(1,1)$. Then 
there are exactly three cylinders of closed trajectories in the
horizontal direction.  
\end{lemma} 

\bold{Proof.}
Each of the two singularities has a saddle connection returning
to itself at angle $\pi$.  Each of these saddles is the boundary
of a cylinder. This accounts for $2$ cylinders. Each boundary
component of a third cylinder consists of the pair of these
saddle connections returning to the singularity. 
\qed\medskip

Note that the width of this wider cylinder is the sum of
the widths of the narrower ones. 
Let $w_1 \|v_1\|,w_2 \|v_1\|$ (where $w_i \in \zed$) be the
widths of
the small cylinders $C_1$ and $C_2$, and let
$w_3 \|v_1\| = (w_1+w_2) \|v_1\|$ be the width of the wide
cylinder
$C_3$.  

We note that the three cylinders may be assembled in two
different
ways.  Either the two narrow cylinders are glued above $z_1$ (and
wide
cylinder is glued below $z_1$), in the sense that as we leave
$z_1$ in
the positive vertical direction we enter one of the narrow
cylinders.
The second possibility is that the wide cylinder is glued above
$z_1$,
and the narrow cylinders are glued below $z_1$.   
Let $\sigma$ be a parameter
which is $+1$ in the first case, and $-1$ in the second.
%-------------------------picture-flatt----------------------------------------
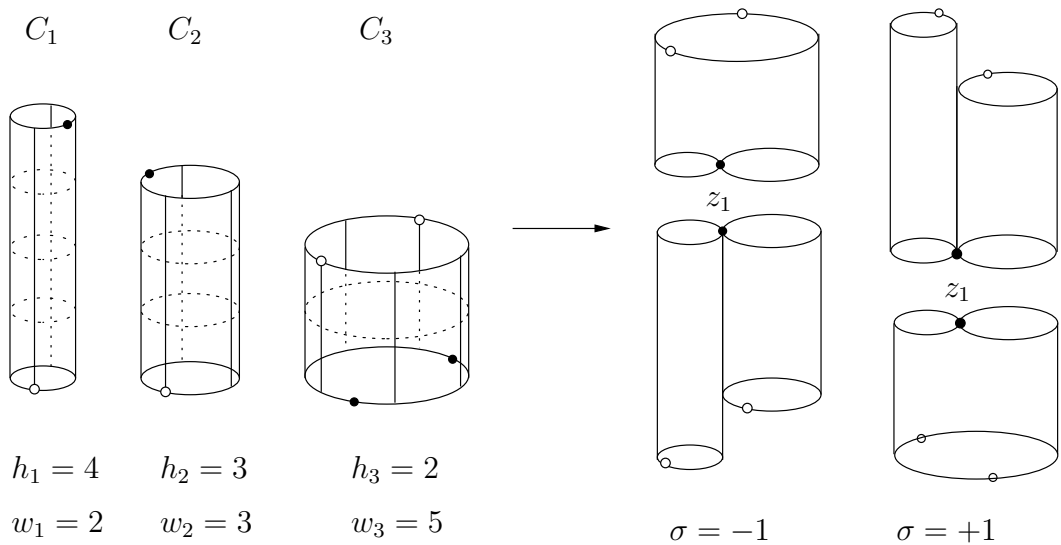
\makefig{How to glue together three cylinders to obtain surfaces of
  genus $2$.}{fig:2}{\begin{picture}(0,0)%
\includegraphics{flatt.pstex}%
\end{picture}%
\setlength{\unitlength}{3947sp}%
\begingroup\makeatletter\ifx\SetFigFont\undefined%
\gdef\SetFigFont#1#2#3#4#5{%
  \reset@font\fontsize{#1}{#2pt}%
  \fontfamily{#3}\fontseries{#4}\fontshape{#5}%
  \selectfont}%
\fi\endgroup%
\begin{picture}(6604,3403)(190,-3269)
\put(6076,-1681){\makebox(0,0)[lb]{\smash{\SetFigFont{12}{14.4}{\rmdefault}{\mddefault}{\updefault}
$z_1$
}}}
\put(218,-2824){\makebox(0,0)[lb]{\smash{\SetFigFont{12}{14.4}{\rmdefault}{\mddefault}{\updefault}
$h_1=4$
}}}
\put(218,-3163){\makebox(0,0)[lb]{\smash{\SetFigFont{12}{14.4}{\rmdefault}{\mddefault}{\updefault}
$w_1=2$
}}}
\put(1151,-2824){\makebox(0,0)[lb]{\smash{\SetFigFont{12}{14.4}{\rmdefault}{\mddefault}{\updefault}
$h_2=3$
}}}
\put(1151,-3163){\makebox(0,0)[lb]{\smash{\SetFigFont{12}{14.4}{\rmdefault}{\mddefault}{\updefault}
$w_2=3$
}}}
\put(2351,-2824){\makebox(0,0)[lb]{\smash{\SetFigFont{12}{14.4}{\rmdefault}{\mddefault}{\updefault}
$h_3=2$
}}}
\put(2351,-3163){\makebox(0,0)[lb]{\smash{\SetFigFont{12}{14.4}{\rmdefault}{\mddefault}{\updefault}
$w_3=5$
}}}
\put(1201,-61){\makebox(0,0)[lb]{\smash{\SetFigFont{12}{14.4}{\rmdefault}{\mddefault}{\updefault}
$C_2$
}}}
\put(2401,-61){\makebox(0,0)[lb]{\smash{\SetFigFont{12}{14.4}{\rmdefault}{\mddefault}{\updefault}
$C_3$
}}}
\put(4576,-1111){\makebox(0,0)[lb]{\smash{\SetFigFont{12}{14.4}{\rmdefault}{\mddefault}{\updefault}
$z_1$
}}}
\put(4351,-3211){\makebox(0,0)[lb]{\smash{\SetFigFont{12}{14.4}{\familydefault}{\mddefault}{\updefault}
$\sigma=-1$
}}}
\put(301,-61){\makebox(0,0)[lb]{\smash{\SetFigFont{12}{14.4}{\rmdefault}{\mddefault}{\updefault}
$C_1$
}}}
\put(5776,-3211){\makebox(0,0)[lb]{\smash{\SetFigFont{12}{14.4}{\familydefault}{\mddefault}{\updefault}
$\sigma=+1$
}}}
\end{picture}
} 
%------------------------------------------------------------------------------
\mc{MARTIN: ARE CHANGES TO FIGURE OK?.} 
Note that each narrow cylinder has a zero at both the top end and
the
bottom end. Let $z'$ denote the intersection with the top end of
the
vertical trajectory passing through the zero at the bottom end.
The
{\em twist} of a narrow cylinder is defined to be the (clockwise)
horizontal distance along the top of the cylinder from $z'$ to
the
zero at the top end.  To define the twist for the wide cylinder
$C_3$,
choose the zero $x$ on the bottom so that as we enter $C_3$
from that zero, the saddle connection on the boundary of $C_1$
lies to the left, and choose the
zero $x'$ on the top so that as we enter $C_3$, the saddle
connection on the boundary of  $C_1$ lies to
the right.  Again let $z'$ denote the intersection with the top
end of
the vertical trajectory passing through $x$.  The {\em twist} is
defined to be the clockwise horizontal distance along the top of
the
cylinder from $z'$ to $x'$. With this definition, if the twist is
$0$, then the
two boundary components of $C_1$ lie directly opposite each other
on
$C_3$.

By construction for $i=1,2$, the twist of the $C_i$ can be
written as
$(t_i + \sigma \delta_1)\|v_1\|$, where $t_i \in \zed$, $0 \le
t_i <
w_i$. The twist of the third (i.e. the wide) cylinder can be
written
as $(t_3 - \sigma \delta_1)\|v_1\|$, where $t_3 \in \zed$, $0 \le
t_3
< w_1+w_2$.

We define the {\em height} of a cylinder to be the length of a
vertical (i.e parallel to $v_2$) trajectory going from the bottom
edge to the top edge.  (This does not represent distance perpenducular
to the horizontal). Similarly, the heights of the cylinder $C_i$ are
of
the form $h_i \|v_2\|$ where the $h_i$ satisfy $h_i \in \zed +
\sigma
\delta_2$ for $i=1, 2$, and $h_3 \in \zed - \sigma \delta_2$.
Clearly
we must also have $h_i >0$. We prefer to use the coordinates $s_1
=
h_1+h_3$, and $s_2 = h_2+h_3$ and $h_3$ instead of $(h_1, h_2,
h_3)$.
Then $1 \le s_1 \in \zed$, $1 \le s_2 \in \zed$, and $0 < h_3 <
\min(s_1, s_2)$. (The upper bound on $h_3$ is determined by the
conditions $h_i = s_i - h_3 >0$, $i=1,2$).

Since the area of $(M,\omega) = 1$, the area of the base torus is
given by $|v_1\times v_2|=1/d$. The area of $C_i$ is $|w_i v_1 \cross
h_i v_2| = h_i w_i |v_1 \cross v_2| = h_i w_i/d$, hence
 we have $w_1 h_1 + w_2 h_2 + w_3 h_3 = d$. 
Using $w_3 = w_1+w_2$ and rewriting we get:
\begin{equation}
\label{eq:area} 
w_1 s_1 +w_2 s_2 =d
\end{equation}

Summarizing this discussion we get the following:
\begin{proposition}
\label{prop:coords:Md:principal}
Let $\Omega$ denote the set of $\{(v_1,v_2), \delta_1, \delta_2,
\sigma, w_1, w_2, s_1, s_2, h_3, t_1, t_2, t_3 \}$ where $(v_1,
v_2)
\in \cF$, $0 \le \delta_1 < 1$, $0 < \delta_2 < 1$, $\sigma \in
\{+1,-1\}$, $w_1$, $w_2$, $s_1$ and $s_2$ are positive integers
satisfying (\ref{eq:area}), $h_3$ is an element of $\zed -\sigma
\delta_2$ satisfying  
$0 < h_3 < \min(s_1, s_2)$, and $t_i$ are integers mod $w_i$,
$i=1,2,3$ (where $w_3 = w_1 + w_2$). Let $\tau: \Omega \to
\Omega$ denote the
map which exchanges  $(w_1, s_1, t_1)$ with $(w_2,
s_2, t_2)$, and fixes the other coordinates.  
Then there is a one to one map $\Psi$ between 
$\Omega/\tau$ and the open dense subset $\cM'_d(1,1)$ of
$\cM_d(1,1)$. The surface $\Psi((v_1,v_2), \delta_1, \delta_2,
\sigma, w_1, w_2, s_1, s_2, h_3, t_1, t_2, t_3)$ corresponding to
a
point of $\Omega$ has the following properties: 
\begin{itemize}
\item It covers the torus $\reals^2/\Delta$ where $\Delta$ is
spanned
  by $(v_1,v_2)$.
\item The point $z_1$ projects to the origin; the point $z_2$
projects 
  to $\delta_1 v_1 + \delta_2 v_2$. 
\item The surface is the union of three cylinders $C_1$, $C_2$
and
  $C_3$. For $i=1, 2$, $C_i$ has width $w_i \|v_1\|$, height
$(s_i -
  h_3)\|v_2\|$ and twist $(t_i + \sigma \delta_1) \|v_1\|$. The
cylinder 
  $C_3$ has width $w_3 = w_1+w_2$, height $h_3$ and twist $(t_3 -
  \sigma \delta_1) \|v_1\|$. 
\item Interchanging $C_1$ and $C_2$ produces the same surface in
  $\cM_d(1,1)$. 
\end{itemize}
\end{proposition}

\begin{corollary}
\label{cor:count:covers}
Let $N_d(1,1)$ denote the number of covers of a fixed torus of
degree
$d$ with $2$ fixed simple branch points. Then, 
\begin{equation}
\label{eq:Nd11}
N_d(1,1) = \sum_{s_1w_1+s_2w_2=d}
w_1w_2(w_1+w_2)\min(s_1,s_2) + \sum_{ 2 s w = d} (2 w^2 s)
\end{equation}
\end{corollary}
\bold{Proof.} We perform the sum over the $t_i$, $h_3$ and
$\sigma$. In the first there is a factor of $2$ from the $2$ choices
of $\sigma$ which cancels the factor of $1/2$ from the action of
$\tau$. However, this does not count correctly the fixed points of
$\tau$, since we should count them with a factor of $1/2$. To
compensate, we add the second term. 

For similar results in a more general setting see
\cite{Dijkgraaf:texel}, \cite{BO}, 
\cite{Eskin:Okounkov:volumes}, and also \cite{Zorich:newton}. 
We note that here we do not weigh each
cover by the inverse of its automorphism group. 
\qed\medskip

The subset $\cH_1(\beta)$ of area $1$ surfaces in  $\cH(\beta)$
carries a 
canonical measure $\nu_1$ defined as follows.
Let $(M,\omega)\in \cH(\beta)$ and let $P$ the set of zeroes of
$\omega$.  \mc{BEGIN CHANGE}
One chooses a basis for the relative homology group
$H_1(M,P,\zed)$. (To avoid confusion, we will henceforth refer to 
the group $H_1(M,\zed)$ as the {\em absolute} homology). \mc{END CHANGE}  
The holonomy of $\omega$ along the relative homology basis gives
local coordinates for $\cH(\beta)$.  One first defines the
measure $\nu$ on $\cH(\beta)$ as the pull-back of  Lebesgue
measure on $\reals^n$ by these local coordinates. 
One checks easily that the measure is well-defined independent of
choice of holonomy basis and is $SL(2,\reals)$ invariant.  Then
for $E\subset \cH_1(\beta)$ define the cone $C(E)$ over $E$ to be
the set  $(M,r\omega)$ where $(M,\omega)\in E$ and $0\leq r\leq
1$. We define 
\begin{displaymath}
\nu_1(E)=n\nu(C(E))
\end{displaymath}

We recall the following: (see \cite{Eskin:Okounkov:volumes},
\cite{Zorich:newton}) 
\begin{proposition}
\label{prop:volume:general}
For any stratum $\cH(\beta)$, 
\begin{displaymath}
\nu_1(\cH_1(\beta)) = n \lim_{D\to\infty}\frac{\sum_{d=1}^D
N_d(\beta)}{D^{n/2}},
\end{displaymath}
where $n = \dim_{\reals} \cH(\beta)$, and $N_d(\beta)$ is the
degree
of the covering map $p: \cM_d(\beta) \to \cT^{|\beta|}$ (i.e. the
number of degree $d$ covers of a fixed torus with ramification
given
by $\beta$). 
\end{proposition}

Combining Corollary~\ref{cor:count:covers} with
Proposition~\ref{prop:volume:general}, we obtain
\begin{lemma}
\label{lemma:volume:H11}
\begin{displaymath}
\nu_1(\cH_1(1,1)) = \frac{\pi^4}{135}.
\end{displaymath}
\end{lemma}

This calculation was first done, using the method presented in
this paper by A.~Zorich. (see \cite{Zorich:newton}).

\bold{Proof of Lemma~\ref{lemma:volume:H11}:} See Appendix~A. 
\qed\medskip

\subsection{The Space  $\cM_d(2)$.}
Now we have a covering map $p: \cM_d(2) \to \cT^1$, where $\cT^1$
is
the space of tori with one marked point, which we always place at
the 
origin. As in \S\ref{sec:Md11}, we may write the base torus as
$\reals^2/\Delta$, where the lattice $\Delta$ is spanned by
$(v_1,
v_2) \in \cF$. Also as in \S\ref{sec:Md11}, we consider the
inverse
image of the horizontal trajectory through the marked point at
the
origin; this breaks up the surface into cylinders. 

Either there are one or two cylinders. We first consider the $2$
cylinder case.  There are a pair of homotopic saddle
connections joining the zero to itself, bounding a cylinder $C_1$
of width $w_1 \|v_1\|$.  There is another saddle connection
returning to
the zero with angle $3\pi$ with length $w \|v_1\|$. Each side of
it
together with one of the previous saddles forms the boundary of a
cylinder $C_2$ which then has width $w_2 \|v_1\|=(w_1+w)\|v_1\|$.

The height of $C_i$ is $h_i \|v_2\|$ where $0 < h_i \in \zed$. 
Also $h_1,h_2$ satisfy 
\begin{equation}
\label{eq:area:2cyl}
h_1 w_1+h_2 w_2=d.  
\end{equation}
The twists $t_i \in \zed$ satisfy $0 \leq t_i<w_i$. 
Hence
\begin{proposition}
\label{prop:H2:2cyl}
Let $\Omega_2 = \{((v_1,v_2), w_1, w_2, h_1,h_2, t_1, t_2) \}$
where 
$(v_1, v_2) \in \cF$, $w_i$ and $h_i$ are positive integers
satisfying 
(\ref{eq:area:2cyl}), $w_1 < w_2$, and $t_i$ are integers
satisfying
$0 < t_i < w_i$. Then there is a bijection $\Psi_2$ between the
elements of $\Omega_2$ and the surfaces of $\cM_d(2)$ which have
$2$
cylinders. The surface $\Psi_2((v_1,v_2), w_1, w_2, h_1,h_2, t_1,
t_2)$ covers the torus $\reals^2/\Delta$ where $\Delta$ is the
lattice 
spanned by $v_1$ and $v_2$. For $i=1,2$, the cylinder $C_i$ has
width $w_i \|v_1\|$, height $h_i \|v_2\|$ and twist $t_i
\|v_1\|$.
\end{proposition}

Now we consider the covers with one cylinder. There are
three closed curves on each
boundary component of the cylinder with lengths $l_1 \|v_1\|,l_2
\|v_1\|,l_3\|v_1\|$ 
such that $l_1+l_2+l_3=w$ and 
$hw=d$, where $h$ is the height. Cyclically permuting the $l_i$
does
not change the surface. Hence,
\begin{proposition}
\label{prop:H2:1cyl}
Let $\Omega_1 = \{((v_1, v_2), l_1, l_2, l_3, h) \}$ where $(v_1,
v_2) 
\in \cF$, $l_1, l_2, l_3, h$ are positive integers satisfying 
$(l_1 + l_2 + l_3) h = d$. Let $\tau: \Omega_1 \to \Omega_1$
denote
the map which cyclically permutes the $l_i$. Then there is a
bijection 
$\Psi_1$ between $\Omega_1/\tau$ and 
the points of $\cM_d(2)$ with one cylinder. The surface
$\Psi_1((v_1,
v_2), l_1, l_2, l_3, h)$ covers the torus $\reals^2/\Delta$ where
$\Delta$ is the lattice  spanned by $v_1$ and $v_2$. The single
cylinder has height $h\|v_2\|$ and width $(l_1+l_2+l_3)\|v_1\|$. 
\end{proposition}

\begin{corollary}
\label{cor:count:H2}
Let $N_d(2)$ denote the number of covers of a fixed torus with
one
double branch point. Then, 
\begin{equation}
\label{eq:count:H2}
N_d(2) = \sum_{\stackrel{h_1w_1+h_2w_2=d}{w_1 < w_2}} w_1w_2 + 
\frac{1}{3}
\sum_{h|d}\sum_{l_1+l_2+l_3=\frac{d}{h}}\frac{d}{h} +
\frac{2}{3}\sum_{3 l = \frac{d}{h}} \frac{d}{h}
\end{equation}
 \end{corollary}
The first term comes from the $2$-cylinder surfaces, the second
from
the $1$-cylinder surfaces with no symmetry, and the third from 
the $1$-cylinder surfaces with $l_1 = l_2 = l_3$ (which have an extra $\zed_3$
symmetry). 

\begin{lemma}
\label{lemma:volume:H2}
\begin{displaymath}
\nu_1(\cH_1(2)) = \frac{\pi^4}{120}
\end{displaymath}
\end{lemma}
\bold{Proof.} This is done by computing the asymptotics of
(\ref{eq:count:H2}) as $d\to \infty$ and using
Proposition~\ref{prop:volume:general}. This again was done by
A.~Zorich, see \cite{Zorich:newton}.
The details of the calculation are in Appendix~A.
\qed\medskip

\subsection{The constants $s_i(1,1)$ and $c(1,1)$.}
\label{sec:generic:constants}
The results of \cite{Eskin:Masur:Zorich:SV} say that  
\begin{displaymath}
s_1(1,1)=3\frac{\nu_1(\cH_1(2))}{\nu_1(\cH_1(1,1))}=
3\frac{\pi^4/120}{\pi^4/135}=\frac{27}{8}, 
\end{displaymath}
and
\begin{displaymath}
s_2(1,1) =
\frac{1}{24}\frac{\nu_1(\cH_1(\emptyset))\nu_1(\cH_1(\emptyset))
}{\nu_1(\cH_1(1,1))} = \frac{5}{8}
\end{displaymath}
The results of \cite{Eskin:Masur:Zorich:SV} also say that 
\begin{displaymath}
c(1,1)=2
\zeta(2)\frac{1}{3}
\frac{\nu_1(\cH_1(\emptyset))}{\nu_1(\cH_1(1,1))}=
2 \zeta(2)\frac{1}{3}\frac{\pi^2/3}{\pi^4/135}= 5
\end{displaymath}
(The number in $\cite{Eskin:Masur:Zorich:SV}$ refer to cylinders
of primitive
geodesics; we need here to multiply by $\zeta(2)$ to count
cylinders of
imprimitive geodesics, and then again by $2$ since when we count
imprimitive geodesics, we count the cylinder in both directions).

\section{Primitive Covers}
\label{sec:primitive}
We say that a cover $\pi: M \to \gt$ is {\em primitive} if it 
does not factor through any other torus cover. 
Clearly to compute the asymptotics of $N_c(S,T)$ and
$N_s^{(i)}(S,T)$  
we may assume that $\pi: M \to \gt$ is a primitive cover. 

We recall the following general fact: 
\begin{lemma}
\label{lemma:primitive:if:basis}
A cover of degree $d$ of a surface of area $d$ in
$\cH(1,1)$ over the standard torus is
primitive if and only if the absolute homology
generates the lattice $\zed\oplus\zed$.
\end{lemma}

\begin{lemma}
\label{lemma:primitive:H11}
Using the coordinates of
Proposition~\ref{prop:coords:Md:principal}, 
the surface \break
$\Psi((v_1,v_2), \delta_1, \delta_2, \sigma, w_1, w_2,
s_1, s_2, h_3, t_1, t_2, t_3)$ is a primitive cover of the torus
$\reals^2/(\zed v_1 + \zed v_2)$ if and only if:
\begin{displaymath}
(s_1,s_2)=1
\end{displaymath}
and
\begin{displaymath}
(s_1 (t_2+t_3)- s_2 (t_1+t_3),w_1,w_2)=1
\end{displaymath}
\end{lemma}

\bold{Proof.} Without loss of generality, we may assume that $v_1
=
(1,0)$, $v_2 = (0,1)$. 
Let $C_1$, $C_2$, $C_3$ be as in
Proposition~\ref{prop:coords:Md:principal}. Recall the
distinguished zeroes $x,x'$ on the bottom and top of $C_3$ used
to define the twist $t_3$. We may take a basis
for the absolute homology as follows:  
\begin{itemize}
\item $a_1$ is the core horizontal curve of  $C_1$.
\item $a_2$  is the core  horizontal curve of $C_2$.
\item $b_1$ is the curve transverse to the horizontal foliation
which
  starts at the zero on the bottom of $C_1$, crosses $C_1$ to the
zero
  on the top and then crosses $C_3$ from $x$ on the bottom to
$x'$ on
  the top.
\item $b_2$ is the curve transverse to the horizontal foliation
which
  starts at the zero on the bottom of $C_2$, crosses $C_2$ to the
zero
  on the top and then crosses $C_3$ from $x$ on the bottom to
$x'$ on
  the top.
\end{itemize}

This absolute homology basis determines holonomy vectors:
\begin{displaymath}
(w_1,0),(w_2,0), (t_1+t_3,s_1), (t_2+t_3,s_2),
\end{displaymath}
(recall that $s_1 = h_1 + h_3$, $s_2 = h_2 + h_3$). 

Let $v= s_1 (t_2+t_3) - s_2 (t_1+t_3)$.
In order to produce a lattice element with second coordinate of
$1$ 
we must have 
$(s_1,s_2)=1$. Let $r=\gcd (w_1,w_2)$.  Every lattice
element of the form $(d,0)$ in the subgroup generated by the 
third and fourth element  satisfies $v|d$. Thus in order to
generate $(1,0)$ it is necessary and sufficient that $(v,r)=1$.
\qed\medskip

\mc{BEGIN CHANGE}

\begin{corollary}
\label{cor:our:surface:is:primitive}
The surface $S(p/q,\alpha)$ is a primitive cover. 
\end{corollary}

\mc{END CHANGE}

\subsection{Connectedness of the space of primitive covers}
Let $\cP_d(\beta) \subset \cM_d(\beta)$ denote the
primitive torus covers in $\cM_d(\beta)$. 
The main result of the section is the following:
\begin{theorem}
\label{theorem:connected}
$\cP_d(1,1)$ is connected. 
\end{theorem}

Theorem~\ref{theorem:asymp} now follows from
Theorem~\ref{theorem:general:asymp},
Corollary~\ref{cor:our:surface:is:primitive} and 
Theorem~\ref{theorem:connected}. 

\bold{Proof of Theorem~\ref{theorem:connected}.}
In this section we will assume that the surface has area $d$, and
the
base torus has area $1$.  
First note that primitivity is 
invariant under continuous deformations. The condition of
primitivity of 
degree $d$ is open: if $(M_n,\omega_n)\to (M,\omega)$ is a
sequence of  degree $d$ covers converging to a primitive cover,
then if $(M_n,\omega_n)$ factored 
through a cover of smaller degree, then so would $(M,\omega)$.
The condition of primitivity is also closed, since the conditions
of Lemma~\ref{lemma:primitive:H11} are clearly closed conditions.
We use the coordinates of
Proposition~\ref{prop:coords:Md:principal}. For any point $S_1$
in
$\cP_d(1,1)$, We can clearly continuously deform the base lattice
$(v_1, v_2)$ until $v_1 = (1,0)$, $v_2 = (0,1)$, i.e. the base
torus
is the standard torus. We can also continuously deform
$(\delta_1,\delta_2)$ to some fixed point $(\epsilon,\epsilon)$.
Thus
any surface is connected by a continuous path to a cover $\pi: S
\to
\bbT^2$, such that $\pi(z_1)=0$, $\pi(z_2) =
(\epsilon,\epsilon)$.
Recall that we have a covering map $p: \cP_d(1,1) \to \cT^2$,
where
$\cT^2 \isom SL(2,\reals) \semidirect \reals^2/SL(2,\zed)
\semidirect
\zed^2$ is the space of tori with two marked points, one of which
is
always at the origin. Hence any point in $\cP_d(1,1)$ can be
connected
by a continuous path to a point of the
$p^{-1}(\bbT^2,(\epsilon,\epsilon))$ (i.e. the fiber of $p$ above
the
point $(\bbT^2,(\epsilon,\epsilon)) \in \cT^2$). We will now show
that any two points in $p^{-1}(\bbT^2, (\epsilon,\epsilon))$ can
be
connected by a continuous path in $\cP_d(1,1)$.  Note that the
points
of $p^{-1}(\bbT^2, (\epsilon,\epsilon))$ are parametrized by the 
discrete parameters $(w_1,w_2, s_1, s_2, h_3, t_1, t_2, t_3,
\sigma)$ of
Proposition~\ref{prop:coords:Md:principal}. 

We now consider the
following three kinds of continuous paths connecting points in
$p^{-1}(\bbT^2,(\epsilon,\epsilon))$:

\bold{The horizontal kernel foliation.}
Let $\gamma_h: \cT^2 \to \cT^2$ be the path
$\gamma_h(t) = (\bbT^2,(\epsilon+t, \epsilon))$, so that
$\gamma_h(0) = \gamma_h(1) = (\bbT^2,(\epsilon, \epsilon))$. 
Since the covering $p$ is unbranched away from the points in
$\cT^2$
where the two marked points coincide, $\gamma_h$ lifts to a 
path $\tilde{\gamma_h}$ on $\cP_d(1,1)$ connecting two points of
the fiber 
$p^{-1}(\bbT^2,(\epsilon,\epsilon))$. Let $F_h:
p^{-1}(\bbT^2,(\epsilon,\epsilon)) \to
p^{-1}(\bbT^2,(\epsilon,\epsilon))$
denote the map taking $\tilde{\gamma_h}(0)$ to
$\tilde{\gamma_h}(1)$. 
It is clear that $F_h$ preserves the leaves of the horizontal
foliations, hence it preserves the heights and widths of the
cylinders, i.e $F_h$ fixes the parameters $w_1, w_2, s_1, s_2,
h_3, \sigma$. 
The twists are
changed by $\pm ||v_1||$ (and thus the $t_i$ are changed by $\pm
1$),  
with the sign determined by $\sigma$. More
precisely, $F_h(\dots, t_1, t_2, t_3, \sigma) = (\dots,
t_1+\sigma, t_2 +
\sigma, t_3 - \sigma, \sigma)$. 

\bold{The vertical kernel foliation.}
Let $\gamma_v: \cT^2 \to \cT^2$ be the path
$\gamma_v(t) = (\bbT^2,(\epsilon, \epsilon + t))$, so that
$\gamma_v(0) = \gamma_v(1) = (\bbT^2,(\epsilon, \epsilon))$. 
Let $\tilde{\gamma_v}$ denote any lift of $\gamma_v$ to
$\cP_d(1,1)$, 
and let $F_v: p^{-1}(\bbT^2,(\epsilon,\epsilon)) \to
p^{-1}(\bbT^2,(\epsilon,\epsilon))$ denote the map sending
$\tilde{\gamma_v}(0)$ to 
$\tilde{\gamma_v}(1)$.

Without loss of generality we assume $w_1 < w_2$. 
Suppose first that $\sigma = +1$. Then as long as $h_3 > 1$, 
$F_v$ sends $h_3$ to $h_3-1$, 
sends $h_i$ to $h_i+1,i=1,2$ 
and fixes all the other coordinates. In
other words the heights of two narrow cylinders increase by $1$,
and the height of the wide cylinder decreases by $1$. If $h_3 <
1$
then we completely collapse the wide cylinder, and  
the cylinder decomposition may change. 
There are in fact two possibilities, depending on the values 
of the twist parameters. If
(when the wide cylinder has height $0$) the top edge of $C_1$
touches the bottom edge of $C_1$, (see Figure~\ref{fig:vmove1}) 
then the new cylinders
have width $(w_1, w_2, w_1+w_2)$ as before, except that we
switched
$\sigma$ to $-\sigma$, and $h_3$ becomes $\epsilon$.
Also $h_1$ goes from $n+\epsilon$ to $n+1-\epsilon$ for some
integer $n$ and similarly for $h_2$. 

%-----------------------------------------------------------------------------------------
\makefignocenter{
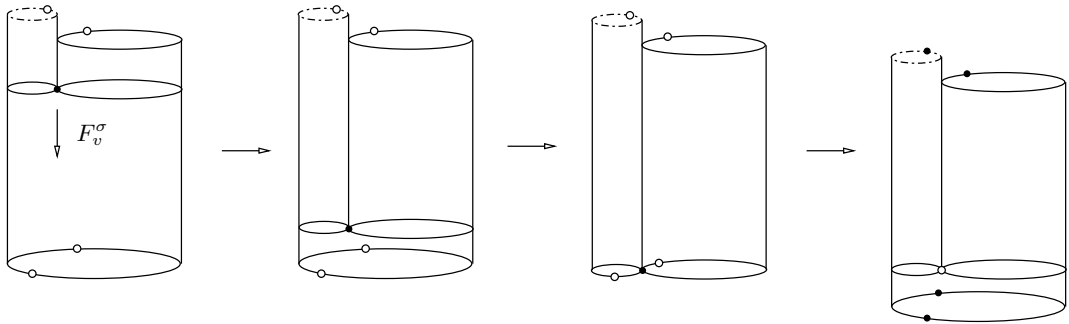
Applying repeatedly $F^{\sigma}_v$ does not change the
cylinder decomposition if $t_3-\sigma\delta_1<w_1$, 
but $\sigma$ changes sign.
}{fig:vmove1}{\begin{picture}(0,0)%
\epsfig{file=verticalmove1.pstex}%
\end{picture}%
\setlength{\unitlength}{4144sp}%
\begingroup\makeatletter\ifx\SetFigFont\undefined%
\gdef\SetFigFont#1#2#3#4#5{%
  \reset@font\fontsize{#1}{#2pt}%
  \fontfamily{#3}\fontseries{#4}\fontshape{#5}%
  \selectfont}%
\fi\endgroup%
\begin{picture}(6355,1904)(259,-2165)
\put(688,-1067){\makebox(0,0)[lb]{\smash{\SetFigFont{9}{10.8}{\rmdefault}{\mddefault}{\updefault}
$F^{\sigma}_v$
}}}
\end{picture}
}
%----------------------------------------------------------------------------------------- 

Suppose on the other hand the
bottom edge of $C_1$ is disjoint from the top edge of
$C_1$.  this corresponds to the twist satisfying $t_3-
\sigma\delta_1>w_1$. 
Then at the stage at which $h_3=0$, the
surface decomposes into cylinders $C_1$ and $C_2$. Each  boundary
component of $C_1$ is the single saddle returning to a zero,
while each boundary component of $C_2$ consists of one of these
saddles as well as a pair of saddle connections joining the two
zeroes. Then as we continue the deformation, $C_2$ becomes the
wide cylinder
and the cylinders have widths $(w_2-w_1,w_1,w_2)$, (see
Figure~\ref{fig:vmove2}). In this case $\sigma$ stays $+1$ and 
if original narrow cylinders had heights
$h_1=n_1+\epsilon$ and $h_2=n_2+
\epsilon$  and the wide cylinder had height $h_3= 1-\epsilon$, 
then after deformation, the heights are $h_1=n_1+1+\epsilon$,
$h_2=\epsilon$ and  $h_3=n_2+1-\epsilon$. 
%-----------------------------------------------------------------------------------------
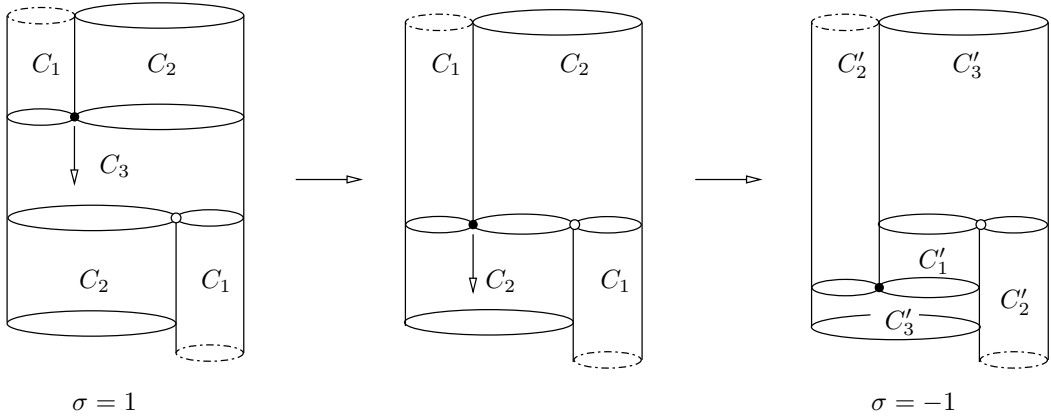
\makefig{Applying repeatedly $F^{\sigma}_v$ changes the cylinder
  decomposition\\ 
and $\sigma$, if $t_3-\sigma\delta_1>w_1$.}{fig:vmove2}{
\begin{picture}(0,0)%
\epsfig{file=verticalmove2.pstex}%
\end{picture}%
\setlength{\unitlength}{4144sp}%
\begingroup\makeatletter\ifx\SetFigFont\undefined%
\gdef\SetFigFont#1#2#3#4#5{%
  \reset@font\fontsize{#1}{#2pt}%
  \fontfamily{#3}\fontseries{#4}\fontshape{#5}%
  \selectfont}%
\fi\endgroup%
\begin{picture}(6253,2486)(-221,-1810)
\put(-62,258){\makebox(0,0)[lb]{\smash{\SetFigFont{10}{12.0}{\rmdefault}{\mddefault}{\updefault}
$C_1$
}}}
\put(2652,-1036){\makebox(0,0)[lb]{\smash{\SetFigFont{10}{12.0}{\rmdefault}{\mddefault}{\updefault}
$C_2$
}}}
\put(181,-1765){\makebox(0,0)[lb]{\smash{\SetFigFont{10}{12.0}{\rmdefault}{\mddefault}{\updefault}
$\sigma=1$
}}}
\put(4960,-1765){\makebox(0,0)[lb]{\smash{\SetFigFont{10}{12.0}{\rmdefault}{\mddefault}{\updefault}
$\sigma=-1$
}}}
\put(3097,258){\makebox(0,0)[lb]{\smash{\SetFigFont{10}{12.0}{\rmdefault}{\mddefault}{\updefault}
$C_2$
}}}
\put(5730,-1158){\makebox(0,0)[lb]{\smash{\SetFigFont{10}{12.0}{\rmdefault}{\mddefault}{\updefault}
$C'_2$
}}}
\put(5243,-916){\makebox(0,0)[lb]{\smash{\SetFigFont{10}{12.0}{\rmdefault}{\mddefault}{\updefault}
$C'_1$
}}}
\put(5041,-1280){\makebox(0,0)[lb]{\smash{\SetFigFont{10}{12.0}{\rmdefault}{\mddefault}{\updefault}
$C'_3$
}}}
\put(4758,258){\makebox(0,0)[lb]{\smash{\SetFigFont{10}{12.0}{\rmdefault}{\mddefault}{\updefault}
$C'_2$
}}}
\put(5446,258){\makebox(0,0)[lb]{\smash{\SetFigFont{10}{12.0}{\rmdefault}{\mddefault}{\updefault}
$C'_3$
}}}
\put(2328,258){\makebox(0,0)[lb]{\smash{\SetFigFont{10}{12.0}{\rmdefault}{\mddefault}{\updefault}
$C_1$
}}}
\put(3340,-1036){\makebox(0,0)[lb]{\smash{\SetFigFont{10}{12.0}{\rmdefault}{\mddefault}{\updefault}
$C_1$
}}}
\put(627,258){\makebox(0,0)[lb]{\smash{\SetFigFont{10}{12.0}{\rmdefault}{\mddefault}{\updefault}
$C_2$
}}}
\put(343,-349){\makebox(0,0)[lb]{\smash{\SetFigFont{10}{12.0}{\rmdefault}{\mddefault}{\updefault}
$C_3$
}}}
\put(951,-1036){\makebox(0,0)[lb]{\smash{\SetFigFont{10}{12.0}{\rmdefault}{\mddefault}{\updefault}
$C_1$
}}}
\put(222,-1036){\makebox(0,0)[lb]{\smash{\SetFigFont{10}{12.0}{\rmdefault}{\mddefault}{\updefault}
$C_2$
}}}
\end{picture}
}
%\vspace*{4mm}
%\begin{center}
%Figure ?: Applying repeatedly $F^{\sigma}_v$ changes the cylinder decomposition\\
%and $\sigma$, if $t_3-\sigma\delta_1>w_1$. 
%----------------------------------------------------------------------------------------- 

In the case $\sigma = -1$, the same analysis applies with $F_v$
replaced by $F_v^{-1}$.

\bold{The $SL(2,\zed)$ action.}
Given an element $g \in SL(2,\zed)$ we now define a map $F_g: 
p^{-1}(\bbT^2,(\epsilon,\epsilon)) \to
p^{-1}(\bbT^2,(\epsilon,\epsilon))$. There is a slight technical
complication since the action of $g \in SL(2,\zed) \subset
SL(2,\reals)$ preserves $\bbT^2$ but not $(\epsilon,\epsilon) \in
\bbT^2$. We solve this as follows: let $(\epsilon',
\epsilon')$ be a point in $\reals^2$ such that $\|g (\epsilon',
\epsilon')\| < 1/3$; then let $\gamma_1(t) = (\bbT^2,
(1-t)(\epsilon, 
\epsilon)+t(\epsilon',\epsilon'))$, $\gamma_2(t) = g(t) ( \bbT,
(\epsilon', \epsilon'))$ where $g(t) \in SL(2,\reals)$ is a path
connecting the identity to $g \in SL(2,\zed)$, and $\gamma_3(t) =
(\bbT^2, (1-t)g(\epsilon', \epsilon')+t(\epsilon,\epsilon))$. Then
let $\gamma$ be the composition of $\gamma_1$,  $\gamma_2$ and
$\gamma_3$; this is a closed path on
$\cT^2$. As above, let $\tilde{\gamma}$ be any lift of $\gamma$,
and
let $F_g: p^{-1}(\bbT^2,(\epsilon,\epsilon)) \to
p^{-1}(\bbT^2,(\epsilon,\epsilon))$ denote the map sending
$\tilde{\gamma}(0)$ to $\tilde{\gamma}(1)$. (A priori the map $F_g$
may depend on some choices; however it will be independent of the
choices in the case where we use it). 
\medskip

We now continue the proof of Theorem~\ref{theorem:connected}.
If $w_1 < w_2$, then 
by applying a suitable power of $F_h$ we may achieve the
situation 
that if the height of the wide cylinder  is collapsed to $0$, the
top and
bottom edge of $C_1$ are disjoint. 
Specifically,  we apply $F_h$ so that $t_3-\sigma\delta_1>w_1$.  
Then we apply $F_v^\sigma$
until the
cylinder decomposition changes. We obtain a surface with three
narrower cylinders. We can now repeat this procedure until $w_1 =
w_2$. 
Thus after doing the moves repeatedly one can assume that the
given torus covering has two 
cylinders of equal width $w$. 

Assuming this, we proceed by applying a power of $F_h$ 
so that 
any vertical trajectory moves from the first cylinder into the
second one (the cylinders are of equal width now) 
and the other way round after returning to the third cylinder.   
After developing the surface in the plane as in Figure~\ref{fig:mov4}, 
we obtain a parallelogram $R$ 
of width $w$ and area $d$ (see the right of Figure~\ref{fig:mov4}) 
and opposite sides identified, 
together with two slits whose endpoints are the singular points.  
That is, the new surface 
consists of a torus $\reals^2/L$ with 2 slits, 
(where $L \subset \zed \oplus \zed$ is a lattice).  We can assume the
base
of one slit is at the origin and the base of the other at an
integer point. After applying a suitable element of the vertical kernel
foliation, we can arrange for both slits have 
length $\epsilon\sqrt{2}$. 
%-------------------------developing the surface-----------------------------------------
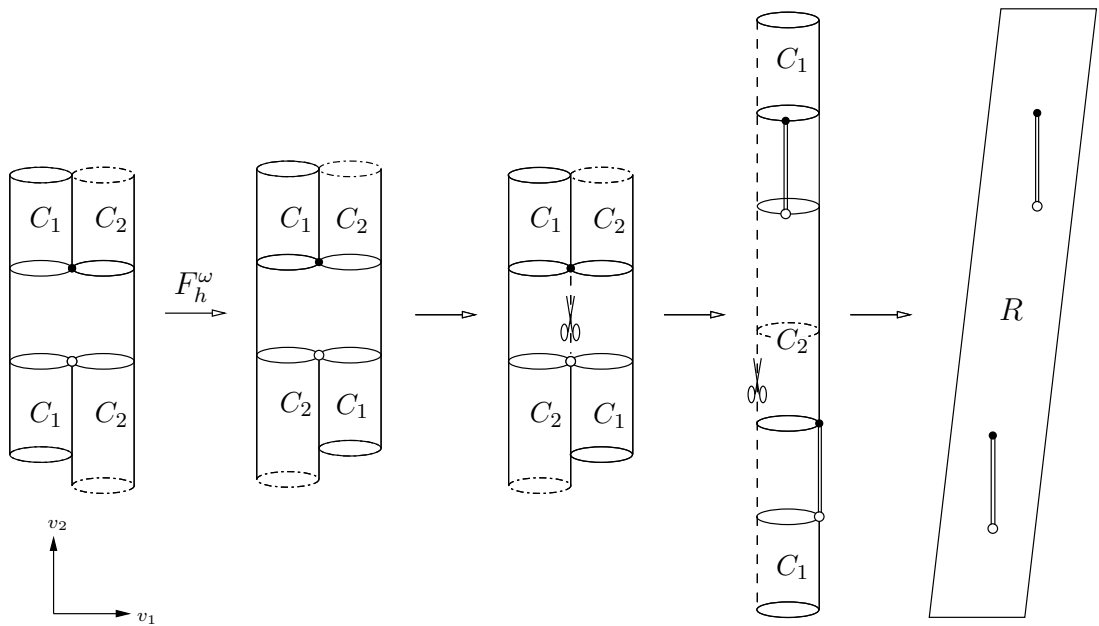
\makefignocenter{Deform and develop a given surface with two narrow cylinders
  of equal width. The twist in the torus on the right depends on the
  twists of the cylinders $C_1$ and $C_2$ (not shown on the left).}{fig:mov4}{\begin{picture}(0,0)%
\includegraphics{move4.pstex}%
\end{picture}%
\setlength{\unitlength}{3947sp}%
\begingroup\makeatletter\ifx\SetFigFont\undefined%
\gdef\SetFigFont#1#2#3#4#5{%
  \reset@font\fontsize{#1}{#2pt}%
  \fontfamily{#3}\fontseries{#4}\fontshape{#5}%
  \selectfont}%
\fi\endgroup%
\begin{picture}(7224,3858)(364,-3157)
\put(512,-678){\makebox(0,0)[lb]{\smash{\SetFigFont{11}{13.2}{\rmdefault}{\mddefault}{\updefault}
$C_1$
}}}
\put(1420,-1097){\makebox(0,0)[lb]{\smash{\SetFigFont{12}{14.4}{\rmdefault}{\mddefault}{\updefault}
$F^{\omega}_{h}$
}}}
\put(1186,-3157){\makebox(0,0)[lb]{\smash{\SetFigFont{7}{8.4}{\rmdefault}{\mddefault}{\updefault}
$v_1$
}}}
\put(626,-2568){\makebox(0,0)[lb]{\smash{\SetFigFont{7}{8.4}{\rmdefault}{\mddefault}{\updefault}
$v_2$
}}}
\put(5202,-2863){\makebox(0,0)[lb]{\smash{\SetFigFont{11}{13.2}{\rmdefault}{\mddefault}{\updefault}
$C_1$
}}}
\put(5202,321){\makebox(0,0)[lb]{\smash{\SetFigFont{11}{13.2}{\rmdefault}{\mddefault}{\updefault}
$C_1$
}}}
\put(4054,-1911){\makebox(0,0)[lb]{\smash{\SetFigFont{11}{13.2}{\rmdefault}{\mddefault}{\updefault}
$C_1$
}}}
\put(3655,-1907){\makebox(0,0)[lb]{\smash{\SetFigFont{11}{13.2}{\rmdefault}{\mddefault}{\updefault}
$C_2$
}}}
\put(4054,-684){\makebox(0,0)[lb]{\smash{\SetFigFont{11}{13.2}{\rmdefault}{\mddefault}{\updefault}
$C_2$
}}}
\put(3655,-678){\makebox(0,0)[lb]{\smash{\SetFigFont{11}{13.2}{\rmdefault}{\mddefault}{\updefault}
$C_1$
}}}
\put(2436,-684){\makebox(0,0)[lb]{\smash{\SetFigFont{11}{13.2}{\rmdefault}{\mddefault}{\updefault}
$C_2$
}}}
\put(2083,-678){\makebox(0,0)[lb]{\smash{\SetFigFont{11}{13.2}{\rmdefault}{\mddefault}{\updefault}
$C_1$
}}}
\put(2436,-1853){\makebox(0,0)[lb]{\smash{\SetFigFont{11}{13.2}{\rmdefault}{\mddefault}{\updefault}
$C_1$
}}}
\put(2083,-1840){\makebox(0,0)[lb]{\smash{\SetFigFont{11}{13.2}{\rmdefault}{\mddefault}{\updefault}
$C_2$
}}}
\put(923,-1907){\makebox(0,0)[lb]{\smash{\SetFigFont{11}{13.2}{\rmdefault}{\mddefault}{\updefault}
$C_2$
}}}
\put(512,-1907){\makebox(0,0)[lb]{\smash{\SetFigFont{11}{13.2}{\rmdefault}{\mddefault}{\updefault}
$C_1$
}}}
\put(923,-678){\makebox(0,0)[lb]{\smash{\SetFigFont{11}{13.2}{\rmdefault}{\mddefault}{\updefault}
$C_2$
}}}
\put(6601,-1261){\makebox(0,0)[lb]{\smash{\SetFigFont{12}{14.4}{\rmdefault}{\mddefault}{\updefault}
$R$
}}}
\put(5202,-1445){\makebox(0,0)[lb]{\smash{\SetFigFont{11}{13.2}{\rmdefault}{\mddefault}{\updefault}
$C_2$
}}}
\end{picture}
}
%----------------------------------------------------------------------------------------- 
%\ref{dev}
 
\begin{lemma}
\label{lemma:double:coset}
Let $P_n$ denote the set of 2 by 2
integral matrices of determinant $n$. Then
representatives for the double cosets
$SL(2,\zed) \backslash P_n/SL(2,\zed)$ consist of the
diagonal matrices
$\begin{pmatrix}
d_1 & 0 \\
0 & d_2 \\
\end{pmatrix}$
where $d_1 d_2 = n$ and $d_1|d_2$. 
\end{lemma}

\bold{Proof.} This is the well known elementary divisor theorem,   
see \cite{Jacobson:rings}. 
\qed\medskip

\begin{corollary}
\label{cor:sl2:orbit}
Let $L$ be a sublattice of $\zed\oplus \zed$ of index $n$. Then
there exists  $g \in SL(2,\zed)$ such that $L' = gL$ is a lattice
such that the fundamental domain for $L'$ is
a rectangle with the length of the short side dividing
the length of the long side and the product of side lengths being
$n$.  
\end{corollary}

\bold{Proof of corollary:} This is just
Lemma~\ref{lemma:double:coset}
in a different language. 
\qed\medskip

We now apply the map $F_g$, 
where $g \in SL(2,\zed)$ is as in Corollary~\ref{cor:sl2:orbit}.
(We note that in this situation, $F_g$ does not
depend on the choice of paths, because the effect of a different path
is to rotate the slit by an angle of $2\pi$.) 
Let
$S'$ denote the resulting surface. Then $S'$ is a rectangular
torus
$\reals^2/L'$, with two slits, one of which is assumed to be at
the
origin. Hence the one slit is  from $(0,0)$ to
$(\epsilon,\epsilon)$, and the
other is from $(u_1,u_2)$ to $(u_1+\epsilon, u_2+\epsilon)$,
where
$u_i \in \zed$. (see Figure~\ref{fig:rectangle}). 
%-------------------------The rectangle---------------------------------------
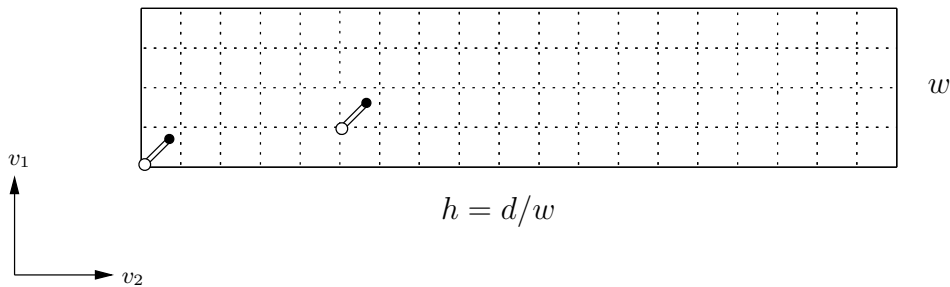
\makefig{The torus $\reals^2/L'$ with two slits developed in
  $\reals^2$.}{fig:rectangle}{\begin{picture}(0,0)%
\epsfig{file=rectangle.pstex}%
\end{picture}%
\setlength{\unitlength}{4144sp}%
\begingroup\makeatletter\ifx\SetFigFont\undefined%
\gdef\SetFigFont#1#2#3#4#5{%
  \reset@font\fontsize{#1}{#2pt}%
  \fontfamily{#3}\fontseries{#4}\fontshape{#5}%
  \selectfont}%
\fi\endgroup%
\begin{picture}(5501,1650)(486,-1602)
\put(3077,-1223){\makebox(0,0)[lb]{\smash{\SetFigFont{12}{14.4}{\rmdefault}{\mddefault}{\updefault}
$h=d/w$
}}}
\put(5987,-464){\makebox(0,0)[lb]{\smash{\SetFigFont{12}{14.4}{\rmdefault}{\mddefault}{\updefault}
$w$
}}}
\put(486,-897){\makebox(0,0)[lb]{\smash{\SetFigFont{9}{10.8}{\rmdefault}{\mddefault}{\updefault}
$v_1$
}}}
\put(1167,-1602){\makebox(0,0)[lb]{\smash{\SetFigFont{9}{10.8}{\rmdefault}{\mddefault}{\updefault}
$v_2$
}}}
\end{picture}
}
%----------------------------------------------------------------------------- 

Since $S'$ was obtained from the original surface $S$ by a
continuous
path, and $S$ was assumed primitive, $S'$ is also primitive. 
Let $u = (u_1,u_2)$ be the vector connecting the 2 slits. Thus
the absolute homology of the surface is generated
by $L'$ and $u$ and the slit curve. The holonomy
along the two sides of the slit curve add to $0$. 
 Thus, by Lemma~\ref{lemma:primitive:if:basis},
$(d_1,0)$, $(0,d_2)$ and $u$ generate $\zed\oplus \zed$.
But $d_1 | d_2$. This is a contradiction unless
$d_1 = 1$. Indeed if $d_1 > 1$, 
reducing modulo $d_1$ we get that
$u \ ({\rm mod}\ d_1)$ generates the finite group $(\zed/d_1
\zed) \oplus
(\zed/d_1 \zed)$:
but that finite group is not generated by any one element. 

So $d_1 = 1$, and hence the surface $S'$ consists of
a vertical  strip of width  $1$, height $d$ 
with $2$ slits, one from $(0,0)$ to $(\epsilon,\epsilon)$ and one
from 
$(0, u_2)$ to $(\epsilon, u_2+\epsilon)$. Since $S'$ is
primitive, 
$(1,0)$, $(0,u_2)$, and $(0,d)$ generate $\zed \oplus \zed$,
hence
$u_2$ and $d$ are relatively prime. Hence there exists $k \in
\zed$
such that $k u_2 = 1\ ({\rm mod}\ d)$.  Now after we apply $F_g$
where 
$g = \begin{pmatrix}
1 & 0 \\
k & 1 \\
\end{pmatrix}$
we get a surface $S_0$ which is again $\reals^2/L'$ with two
slits,
one from $(0,0)$ to $(\epsilon,\epsilon)$ and one from $(0,1)$ to
$(\epsilon,1+\epsilon)$. Hence we have joined any primitive cover
to a 
fixed cover $S_0$ by a continuous path of primitive covers. 
This completes the proof of Theorem~\ref{theorem:connected}.
\qed 

%\ref{dev}

\subsection{Counting primitive covers}
Let  $N_d^P(1,1)$ denote the number of  
of primitive covers of degree $d$ of a surface
of genus $2$ that are branched over $2$ points of the standard
torus and  $N_d^P(2)$,  the number  of
primitive covers of degree $d$ of a surface of genus $2$ branched
over a single point. 

Let $v$ be the quantity 
\begin{displaymath}
s_2(t_1+t_3)-(t_2+t_3)s_1=(t_1+t_3,s_1)\times (t_2+t_3,s_2).
\end{displaymath}

\begin{lemma}
\label{lemma:twists}
Assume $(s_1,s_2)=1$.  For fixed $w_1,w_2$  the number of 
$(t_1,t_2,t_3)$ satisfying  $0 \le t_1 <w_1$, $0 \le t_2 < w_2$,
$0 \le t_3 <
 w_1+w_2$, $r|w_1,r|w_2,r|v$ is 
\begin{displaymath}
\frac{w_1w_2(w_1+w_2)}{r}
\end{displaymath}
\end{lemma}
\bold{Proof.}
Define a linear map 
\begin{displaymath}
L:\zed/r\zed\times
\zed/r\zed\times\zed/r\zed\to\zed/r\zed
\end{displaymath}
by \begin{displaymath}
L(t_1,t_2,t_3)=v\ mod(r)
\end{displaymath}
Since $(s_1,s_2)=1$, the map $L$ is onto. Thus $|ker L|=r^2$. 
Now
dividing the intervals $[0,w_1]$, $[0,w_2],[0,w_1+w_2]$ into
subintervals of length $r$, in each triple of subintervals, we
have
exactly $r^2$ solutions of $r|v$.  Since there are
$w_1w_2(w_1+w_2)/r^3$ triples of intervals, the lemma follows.
\qed\medskip

\begin{lemma}
\label{lemma:NPd:11}
We have $N_2^P(1,1) = 4$. If $d \ge 3$, 
\begin{equation}
\label{eq:NdP:11}
N_d^P(1,1) = \sum_{r|d} \mu(r)
\left(\sum_{\stackrel{(s_1,s_2)=1}{s_1u_1+s_2u_2=\frac{d}{r}}}
r^2u_1u_2(u_1+u_2) \min(s_1,s_2) \right),
\end{equation}
where $\mu(\cdot)$ is the Mobius function. 
\end{lemma}

\bold{Proof.} We use the coordinates of 
Proposition~\ref{prop:coords:Md:principal}, with $v_1 = (1,0)$,
$v_2 = 
(0,1)$.

Note that in the proposition the torus was not standard, while in
counting covers, we can assume that the torus is standard.
Since we are counting primitive covers, we can in view of
Lemma~\ref{lemma:primitive:H11} from now on assume
\begin{displaymath}
(s_1,s_2)=1
\end{displaymath}  
By Lemma~\ref{lemma:twists} the number of covers that satisfy
$r|w_1,r|w_2,r|v$ is 
\begin{displaymath}
\sum_{\stackrel{s_1w_1+s_2w_2=d}{r|w_1,r|w_2}}
\frac{w_1w_2(w_1+w_2)}{r}\min(s_1,s_2) +  \sum_{r | w}
\sum_{2 w = d} 2w^2 
\end{displaymath}
(where the last term is coming from the fixed points of $\tau$.
Indeed 
$(s_1,s_2) = 1$ and $s_1 = s_2$ together imply $s_1 = s_2 = 1$,
and
since $t_1 = t_2$, $v = 0$. 
This means that there is no condition on the $t_i$ other than
$0 \le t_i<w_i$.  
Also $w_1 = w_2 = w$, hence $d = 2 w$. 
Thus the number of fixed points of $\tau$ satisfying $(s_1,s_2) =
1$ and $r|w_1$, $r|w_2$, $r|v$ is 
$2 w^2 = \tfrac{d^2}{2}$ if
$r|\tfrac{d}{2}$ and $0$ otherwise.) 

After substituting $w_i = r u_i$ the sum becomes 
\begin{equation}
\label{eq:H11:fixed:r}
\left(\sum_{\stackrel{(s_1,s_2)=1}{s_1u_1+s_2u_2=\frac{d}{r}}}
r^2u_1u_2(u_1+u_2) \min(s_1,s_2)\right) + \chi(r,d)\frac{d^2}{2}
\end{equation}
where $\chi(r,d) = 1$ if $r|\tfrac{d}{2}$ and $\chi(r,d)=0$
otherwise. 
Now, using the Mobius inversion formula, we get
(\ref{eq:NdP:11}). 
(Note that the fact that the contribution of the fixed points 
of $\tau$ cancels can be seen directly since there are no
primitive
covers which are fixed points of $\tau$ unless $d = 2$.) 
\qed\medskip

\begin{lemma}
\label{lemma:asymp:NdP11}
As $d \to \infty$, 
\begin{displaymath}
N_d^P(1,1) = \frac{d^4}{3}\sum_{r|d}\frac{\mu(r)}{r^2} + o(d^4)
\end{displaymath}
\end{lemma}
\bold{Proof.} See Appendix~A. 
\qed\medskip

\bold{Remark:} In fact it can be derived from the results of 
\cite{Dijkgraaf:texel} and \cite{BO} that
\begin{displaymath}
N_d^P(1,1) = \frac{1}{3}d^3(d-1)\sum_{r|d}\frac{\mu(r)}{r^2}
\end{displaymath}

\begin{lemma}
\label{lemma:NP:2}
Let $N_d^P(2)$ denote the number of primitive genus $2$ covers of
the standard
torus branched over a single point.  Then, $N_2^P(2)
= 0$, $N_3^P(2) 
= 3$ and for $d \ge 4$, 
\begin{equation}
\label{eq:NP:2}
N_d^P(2) = \sum_{r|d} \mu(r) \left(
\sum_{\stackrel{h_1u_1+h_2u_2=d/r}{\stackrel{(h_1,h_2)=1}{u_1
<u_2}}}
r u_1 u_2 + \frac{1}{3} \sum_{u_1+u_2+u_3=d/r} d \right)
\end{equation}
\end{lemma}

\bold{Proof.} We begin with counting the number of primitive
covers with $2$ cylinders in $\cH(2)$, using the coordinates of
Proposition~\ref{prop:H2:2cyl}. Again we are counting covers over
the standard torus. 
For the cover to be
primitive we must
have $(h_1,h_2)=1$, which we now assume.   We must also have
$(w_1,w_2,(t_1,h_1)\times
(t_2,h_2))=1$.                
As in the lemma, the number of twists such that $r|w_1,
r|w_2,
r|(t_1h_2-t_2h_1)$
is \begin{displaymath}
\frac{w_1w_2}{r}.
\end{displaymath}
Thus the number of primitive covers of degree
$d$ with $2$ cylinders is:
\begin{equation}
\label{eq:NP2:2cyl:r:fixed}
\sum_{r|d} \mu(r)
\sum_{\stackrel{h_1u_1+h_2u_2=d/r}{\stackrel{(h_1,h_2)=1}{u_1
<u_2}}}
r u_1 u_2
\end{equation}
Note that there are no covers in degree $2$.

We do the same counting of the number of primitive covers with
one
cylinder, using the coordinates of
Proposition~\ref{prop:H2:1cyl}. 
We must have $h=1$ in order for the cover to be primitive.
As before we make this assumption. We then have
$d=w=l_1+l_2+l_3$,
where the $l_i$ are the lengths of the curves.  For the cover to
be
primitive we must also have $(l_1,l_2,l_3)=1$.  Then arguing as
before, using the symmetry of the $l_i$, the number of primitive
covers is for $d \ge 4$, 
\begin{equation}
\label{eq:NP2:1cyl:r:fixed}
\frac{1}{3}\sum_{r|d} \mu(r) 
 \sum_{u_1+u_2+u_3=d/r} d.
\end{equation}  
(Note that there are no primitive covers with $l_1=l_2=l_3$
unless $d=3$).
Now the lemma follows immediately from
(\ref{eq:NP2:2cyl:r:fixed}) and (\ref{eq:NP2:1cyl:r:fixed}).
\qed\medskip

\begin{lemma}
\label{lemma:asymp:NdP:2}
As $d \to \infty$, 
\begin{displaymath}
N_d^P(2) = \frac{3}{8} d^4 \sum_{r|d}\frac{\mu(r)}{r^2} +
o(d^3) 
\end{displaymath}
\begin{displaymath}
\end{displaymath}
\end{lemma}
\bold{Proof.}  See Appendix~A. 
\qed\medskip

\bold{Remark:} In fact it can be shown that
\begin{displaymath}
N_d^P(2) = \frac{3}{8} d^2(d-2) \sum_{r|d}\frac{\mu(r)}{r^2}
\end{displaymath}

\subsection{The constant $s_1(d)$.}
We normalize the space of tori to have measure $1$. 
For a surface of area $1$ which is a $d$ fold cover of a torus,
the torus has area $\frac{1}{d}$.  Thus 
\begin{equation}
\label{eq:nuP11:to:NdP11}
\nu(P_d(1,1))=\frac{1}{d}N_d^P(1,1)
\end{equation} while
\begin{displaymath}
\nu(P_d(2))=N_d^P(2)
\end{displaymath}

Given a surface in $\cH(1,1)$ and a saddle connection $\gamma$ joining the
two zeroes of small length, if there is no saddle connection
homologous to $\gamma$,  we can deform the surface by
letting the length of the saddle connection go to $0$. The
resulting surface lies in
$\cH(2)$.  We  keep
the absolute homology constant during the deformation.  Conversely, 
given a surface in  $\cH(2)$, and a vector $\gamma$ of length $\epsilon$,
we can break up the double zero into two simple zeroes together
with a saddle connection joining the two zeroes such that the holonomy
of the flat structure along the saddle is $\gamma$. 
Because the total angle at
the double zero is $6\pi$, there are in fact $3$ ways of doing this
(we refer the reader to \cite{Eskin:Masur:Zorich:SV} for details). 
Since this process preserves the  absolute homology,  
the next lemma is immediate:
\begin{lemma}
\label{lemma:breakup}
Given a surface in $P_d(1,1)$ the surface given by collapsing a single
saddle connection lies in $P_d(2)$.  Conversely, given a surface
in $S \in P_d(2)$, and a short vector $\gamma \in \reals^2$, 
there are $3$ surfaces with two simple zeroes joined by a single
saddle connection of holonomy $\gamma$, which have the property that
after collapsing $\gamma$ we get $S$. All three surfaces
lie in $P_d(1,1)$ (because the absolute homology is preserved). 
\end{lemma}

\begin{proposition}
\label{prop:value:s1}
The constant $s_1(d)$ in Theorem~\ref{theorem:asymp} 
can be expressed in terms of the number of $d$-fold
branched 
covers of the torus with prescribed branching. In fact, we have, 
\begin{displaymath}
s_1(d) = 3\frac{ dN_d^P(2)}{N_d^P(1,1)}
\end{displaymath}
where $N_d^P(2)$ is the number of
$d$-fold primitive covers of a surface of genus $2$ over a torus,
branched over a single
point, and $N_d^P(1,1)$ is the number of $d$-fold
primitive covers over the torus with simple order 2 branching at
two distinct
points. 
\end{proposition}

\bold{Proof of Proposition~\ref{prop:value:s1}}
The   Siegel-Veech formula, Lemma~\ref{lemma:breakup}  and  the
connectedness  of $P_d(1,1)$ (Theorem~\ref{theorem:connected})
give
\begin{displaymath}
s_1(d)=3\frac{\nu(P_d(2))}{\nu(P_d(1,1))}
\end{displaymath}
which by (\ref{eq:nuP11:to:NdP11})
is \begin{displaymath} 3
\frac{dN_d^P(2)}{N_d^P(1,1)},
\end{displaymath}
where $N_d^P(2)$ is given by (\ref{eq:NP:2}) and $N_d^P(1,1)$ is
given 
by (\ref{eq:NdP:11}).

\medskip\noindent

\bold{Proof of Theorem~\ref{theorem:constants:converge} for
$s_1(d)$.}
In view of Lemma~\ref{lemma:asymp:NdP:2}, 
Lemma~\ref{lemma:asymp:NdP11},  and
\S\ref{sec:generic:constants}, 
\begin{displaymath}
\lim_{d\to\infty}s_1(d)=\lim_{d\to\infty} 
 3 \frac{d N_d^P(2)}{N_d^P(1,1)}=
 3(\frac{3/8}{1/3})=\frac{27}{8}=s_1(1,1).
\end{displaymath}
\qed

\subsection{The constant $c(d)$.}

\begin{proposition}
\label{prop:value:c:d}
For $d \ge 3$, we have
\begin{equation}
\label{eq:value:c:d}
c(d) = \frac{d}{N_d^P(1,1)} \sum_{r|d} \mu(r)
\sum_{\stackrel{(s_1,s_2)=1}{s_1u_1+s_2u_2=\frac{d}{r}}}
 u_1u_2(u_1+u_2)\min(s_1,s_2) \left(
\frac{1}{u_1^2} + \frac{1}{u_2^2} + 
  \frac{1}{(u_1+u_2)^2} \right),
\end{equation}
where $N_d^P(1,1)$ is given by (\ref{eq:NdP:11}), 
\end{proposition}

\bold{Proof.} We let $f$ be the characteristic function of a disc
of
radius $\epsilon$ in $\reals^2$. The Siegel Veech formula says
that:
\begin{equation}
\label{eq:c:d:immed}
c(d)\pi\epsilon^2=\frac{\int_{\cP_d(1,1)}\hat f \,
  d\nu}{\nu(P_d(1,1))} = \zeta(2) \frac{\int_{\cP_d(1,1)}\tilde f
\,
  d\nu}{\nu(P_d(1,1))}
\end{equation}
Now $\hat f$ counts the number of cylinders of (imprimitive) 
closed geodesics with 
length at most $\epsilon$, and $\tilde{f}$ counts the number of
primitive cylinders of closed geodesics of length at most
$\epsilon$. 
(The relation
$\zeta(2) \int_{\cP_d(1,1)}\tilde f \, d\nu =
\int_{\cP_d(1,1)}\hat f \, d\nu$ was proved in a more general
setting in
\cite{Eskin:Masur:ae}, see the remark following the proof of
Theorem~5.1(b)).
 
We wish to evaluate the numerator in the above expression.  
As above, we use the coordinates of
Proposition~\ref{prop:coords:Md:principal}. Clearly if
$\tilde{f}(S) \ne 
0$, $S$ contains a closed curve of length at most $\epsilon$. 
Since $s_1, s_2$ are positive integers and $\|v_2\|\ge \|v_1\|$,
there
is a uniform lower bound on the length of any curve which is is
not
parallel to $v_1$. Hence, for $\epsilon$ sufficiently small, the
three 
candidates curves of length at most $\epsilon$ are the widths of
the three cylinders. Hence, (assuming $w_1 \le w_2$) for
sufficiently
small $\epsilon$, 
\begin{displaymath}
\tilde{f}(S) = \begin{cases}
0, & \text{ if $w_1 \|v_1\| > \epsilon$} \\
1, & \text{ if $w_2 \|v_1\| > \epsilon > w_1 \|v_1\|$ } \\
2, & \text{ if $(w_1+w_2) \|v_1\| > \epsilon > w_2 \|v_1\|$ } \\
3, & \text{ if $\epsilon > (w_1+w_2) \|v_1\|$ } \\
\end{cases}
\end{displaymath}
Let $\chi: \reals^2 \to \reals$ be defined by:
\begin{displaymath}
\chi(v) = \begin{cases}
0, & \text{ if $w_1 \|v\| > \epsilon \sqrt{d}$} \\
1, & \text{ if $w_2 \|v\| > \epsilon \sqrt{d} > w_1 \|v\|$ } \\
2, & \text{ if $(w_1+w_2) \|v\| > \epsilon \sqrt{d} > w_2 \|v\|$
} \\
3, & \text{ if $\epsilon \sqrt{d} > (w_1+w_2) \|v\|$ } \\
\end{cases}
\end{displaymath}
Then 
\begin{equation}
\label{eq:int:tilde:chi}
\int_{\reals^2} \chi(v) \, dv = \pi \epsilon^2 d \left(
\frac{1}{w_1^2}
  + \frac{1}{w_2^2} + \frac{1}{(w_1+w_2)^2} \right)
\end{equation}
For a unimodular lattice $\Delta \subset \reals^2$, let $\Delta'$
denote the primitive vectors in $\Delta$, and let
$\tilde{\chi}(\Delta) = \sum_{v \in \Delta'} \chi(v)$. Then, by
the
Siegel formula, 
\begin{equation}
\label{eq:siegel:tilde:chi}
\int_{SL(2,\reals)/SL(2,\zed)} \tilde{\chi}(\Delta) \,
d\mu(\Delta) =
\frac{1}{\zeta(2)} \int_{\reals^2} \chi(v) \, dv
\end{equation}
and $\mu$ is normalized Haar measure on the space of unimodular
tori
$SL(2,\reals)/SL(2,\zed)$.
Note that if we fix all the parameters except $(v_1,v_2)$, 
\begin{equation}
\label{eq:tilde:f:is:tilde:chi}
\tilde{f}((v_1,v_2), \dots) = \tilde{\chi}(\zed v_1 \sqrt{d}
\oplus
\zed v_2 \sqrt{d})
\end{equation}
and $\zed v_1 \sqrt{d} \oplus \zed v_2 \sqrt{d}$ is a unimodular
lattice.  
We now compute $\int_{\cP_d(1,1)} \tilde{f}(S) \, d\nu(S)$ by
parametrising $S$ as in
Proposition~\ref{prop:coords:Md:principal} and
Lemma~\ref{lemma:NPd:11}, and computing the integral over
$(v_1,v_2)$
first. Then, in view of (\ref{eq:tilde:f:is:tilde:chi}), 
(\ref{eq:siegel:tilde:chi}), and (\ref{eq:int:tilde:chi}) we
obtain 
\begin{multline}
\label{eq:numerator}
\int_{\cP_d(1,1)} \tilde{f}(S) \, d\nu(S) = \\
= \frac{\pi \epsilon^2 d}{d \zeta(2)}\sum_{r|d} \mu(r)
\sum_{\stackrel{(s_1,s_2)=1}{s_1u_1+s_2u_2=\frac{d}{r}}} 
u_1u_2(u_1+u_2) \min(s_1,s_2) \left( \frac{1}{u_1^2} +
\frac{1}{u_2^2} + 
  \frac{1}{(u_1+u_2)^2} \right)
\end{multline}
where in the sum $w_i = r u_i$, and the factor of $d$ in the
denominator comes from the integral over the $\delta_i$ in
Proposition~\ref{prop:coords:Md:principal}. 
Now the lemma follows from (\ref{eq:numerator}), 
(\ref{eq:c:d:immed}) and (\ref{eq:nuP11:to:NdP11}).  
\qed\medskip

The proof of Theorem~\ref{theorem:constants:converge} for $c(d)$ is postponed
 until Appendix~A. 

\subsection{The constant $s_2(d)$.}
\begin{proposition}
\label{prop:value:s2:d}
Assume $d \ge 3$. We have 
\begin{equation}
\label{eq:value:s2:d}
s_2(d)=\frac{d}{N_d^P(1,1)}\left[\sum_{r|d} \mu(r) \left(
\sum_{\stackrel{h_1u_1+h_2u_2=d/r}{\stackrel{(h_1,h_2)=1}{u_1
<u_2}}} 
r u_1 u_2 \right) +
\sum_{\stackrel{w|d}{w\ne d}} \varphi(d/w)
\sum_{r|w} \frac{\mu(r)}{r}w^2 + \frac{d \varphi(d)}{2}
\right]
\end{equation}
where $N_d^P(1,1)$ is given by (\ref{eq:NdP:11}), and $\varphi$ is the 
Euler function. 
\end{proposition}
\bold{Proof.}
 We let $f$ be the characteristic function of a disc of
radius $\epsilon$ in $\reals^2$. The Siegel Veech formula says
that:
\begin{equation}
\label{eq:c:d:immed2}
s_2(d)\pi\epsilon^2=\frac{\int_{\cP_d(1,1)}\hat f \,
  d\nu}{\nu(P_d(1,1))}
\end{equation}
Now $\hat f$ counts the number of  
saddle connections of multiplicity $2$ of 
length at most $\epsilon$. 
We wish to evaluate the asymptotics of the numerator as $\epsilon
\to
0$. 

We use the coordinates of
Proposition~\ref{prop:coords:Md:principal}, 
except we choose to work with surfaces of area $d$ (which cover 
a torus of area $1$). 
We first show that the contribution to the integral from the part
where the base torus is degenerate has lower order
$O(\epsilon^4)$. For suppose   $||v_1||\leq \epsilon$ so that the
measure of the set of such $v_1$ is $O(\epsilon^2)$. The
saddle connections cannot be horizontal for they would bound cylinders.  
Thus it must happen that  each
crosses a cylinder which means that the distance across a cylinder
must be bounded by $\epsilon$.  Since this is true for two
cylinders, we have another factor of $O(\epsilon^2)$, coming from the measure of the heights of the cylinders, proving the
claim.

 Since $s_i = h_i + h_3$ are integers, there are 
two ways
to produce two saddle connections of length at most $\epsilon$:
either the pair of short saddle connections cross the wide
cylinder $C_3$ ($h_3 \le \epsilon$) or one saddle crosses $C_1$ ($h_1
\le \epsilon$) 
and the other crosses $C_2$ ($h_2 \le \epsilon$). (This is true
because any curve  
wrapping around a cylinder would have length at least $\epsilon$). 
In either case we have $\delta_2 \leq\epsilon$. 

\bold{Case 1:} The pair crosses $C_3$. 
Recall $x$ is the representation of the zero on the bottom so that
$C_1$ lies to the left, and $x'$ is on the top so that $C_1$ lies to
the right. Let $y,y'$ the other representation of the zeroes on the
bottom and top respectively.  There are  two subcases 
here. 
The first subcase is if one saddle connection joins $x$ on the
bottom to $y'$ on the top.  
In order for there to be a  second saddle  connection joining $y$ to $x'$  parallel of the same length to the first, we 
must have $w_1 = w_2$ and $t_3 = w_1$.
By primitivity $(s_1,s_2)=1$. 
Also by primitivity we have $(s_2(t_1+t_3)-s_1(t_2+t_3),w_1)=1$ or
$(s_2t_1-s_1t_2,w_1)=1$.  
Also the projections of the two zeroes to the base torus must be
at most $\epsilon$ apart. There is also a factor of $1/2$ from the
action of $\tau$. 
Then the contribution to the integral for this subcase is
\begin{displaymath}
\frac{\pi \epsilon^2}{2}\sum_{\stackrel{(s_1+s_2)w=d}{(s_1,s_2)=1}}\quad
\sum_{\stackrel{(s_1t_2 -s_2 t_1,w) = 1}{0\le t_i < w}}1 = \frac{\pi
\epsilon^2}{2}\sum_{\stackrel{w|d}{w\ne d}} \varphi(d/w) 
\sum_{r|w} \frac{\mu(r)}{r}w^2
\end{displaymath}
where $\varphi$ is the Euler function.

In the second subcase, one saddle joins $x$ on the bottom to $x'$ on
the top and the  other saddle connection across $C_3$
joins $y$ to $y'$. They are  parallel of the same length. 
The saddle connections are 
homologous and  together  form a dividing curve. 
If we cut the surface along the dividing curve we have a pair of
tori glued together. 
Conversely, given a pair of tori with marked point at the origin on each, 
and a  vector $\gamma$, we  can realize the surface as the tori glued along the slits.   Namely, we take the vector $\gamma$ based
at the origin, slit each torus along $\gamma$ and glue the slit
tori. 
The tori have widths $w_i$, heights $h_i$ and twists $t_i$. 
Since the area must be $d$ we have
\begin{displaymath}
h_1w_1+h_2w_2=d
\end{displaymath}
Moreover the lattices of the two tori together must  
generate a unimodular lattice since the surface itself was a
primitive cover.
Thus we must have a similar pair of conditions we encountered 
in the case of two cylinders in $\cH(2)$; namely, 
we must have $(h_1,h_2)=1$ and $(w_1,w_2,t_1h_2-t_2h_1)=1$.
On the other hand, $w_1$ and $w_2$ are arbitrary, but interchanging
the two tori yields the same surface. Hence there is a term exactly as 
in the two cylinder case $\cH(2)$ with $w_1 < w_2$, and also a term 
when $w_1 = w_2$ which is exactly the term in subcase 1. 
Again using Theorem~\ref{theorem:connected}
we find that the contribution to the integral from this subcase is given by  
\begin{displaymath}
\pi\epsilon^2 \sum_{r|d} \mu(r)
\sum_{\stackrel{h_1u_1+h_2u_2=d/r}{\stackrel{(h_1,h_2)=1}{u_1
<u_2}}}
r u_1 u_2+\frac{\pi \epsilon^2}{2}\sum_{\stackrel{w|d}{w \ne d}}\varphi(d/w)
\sum_{r|w} \frac{\mu(r)}{r}w^2 
\end{displaymath}
%This gives the first term in the sum. 

\bold{Case 2:} One saddle crosses $C_1$, the other crosses $C_2$.
We must have  $t_1 = t_2 = 0$ in order for the saddle to be
short.  
Again since $s_i = h_i + h_3$, and by primitivity $(s_1,s_2) =
1$, we
must
have $s_1 = s_2 = 1$. Also because of primitivity, $(w_1, w_2) =
1$. 
The area $s_1 w_1 + s_2 w_2 = d$, hence $w_1 + w_2 = d$. Also
$t_3$
can be arbitrary.   There is
a factor of $1/2$ due to the action of $\tau$. 
Hence the contribution to the integral is
\begin{displaymath}
\frac{1}{2}\pi \epsilon^2\sum_{\stackrel{w_1 + w_2 = d}{(w_1,w_2) = 1}} (w_1 +
w_2) = \frac{\pi \epsilon^2 d}{2} \varphi(d)
\end{displaymath}
The Proposition follows by adding the contributions for the three
cases. 
\qed\medskip

The proof of Theorem~\ref{theorem:constants:converge} for $s_2(d)$ is postponed
 until Appendix~A. 

\appendixmode

\section{Appendix: Some asymptotic formulas}
%\label{sec:asymp}
\bold{Proof of Lemma~\ref{lemma:volume:H11}.}
It is easy to see that the second term in (\ref{eq:Nd11}) is
$o(d^4)$. Hence, 
\begin{equation}
\label{eq:asymp:Nd:1:1}
N_d(1,1)
=\sum_{r|d}\sum_{\stackrel
{(s_1,s_2)=1}{s_1w_1+s_2w_2=d/r}}r
w_1w_2(w_1+w_2)\min(s_1,s_2) + o(d^4)
\end{equation}

We compute the asymptotics  as $d\to\infty$. 
Suppose $d,r$ satisfy $d/r\leq d^{1/5}$.
Then $w_i,s_i\leq d^{1/5}$ and so each term is at most $2dd^{4/5}$.  There are at most $dd^{4/5}$ terms so the
contribution to the sum is at most $2d^2d^{8/5}$.  
The
contribution to the sum from just the term  $r=1$ is of the order
$d^4$, as we shall
see.  Thus the asymptotics of the 
contribution from  those $r$ for which 
$r\leq d^{4/5}$
 is
the same as the total contribution. 
In particular $d/r\to\infty$ as $d\to\infty$. 
For these terms, we will calculate the contribution, each with an
error that is $o(d/r)^4)$ so that the total error is $o(d^4)$ as
$d\to\infty$. 

For each  $(s_1,s_2)$, choose
a smallest $w_2^0\leq s_1$ such  that 
\begin{displaymath}
s_1w_1+s_2w_2=d/r
\end{displaymath}
has a solution $(w_1^0,w_2^0)$.
Then for any pair $(w_1,w_2)$ such that $w_1s_1+w_2s_2=d/r$,
subtracting $(w_1^0,w_2^0)$,  
we find 
\begin{displaymath}
s_1(w_1-w_1^0)+s_2(w_2-w_2^0)=0
\end{displaymath}
 and since
$(s_1,s_2)=1$, there must be $\lambda$ 
so that  
\begin{displaymath}
(w_1,w_2)=(w_1^0,w_2^0)+\lambda (-s_2,s_1)
\end{displaymath}
with  
\begin{equation}
\label{eq:lambda:range}
0\leq\lambda\leq \frac{d/r}{s_1s_2}
\end{equation}

Note that $(w_1 w_2)(s_1 s_2) = (w_1 s_1 )(w_2 s_2) \le
(d/r)^2$. Hence, 
\begin{equation}
\label{eq:trivial:bound}
w_1 w_2 \le \frac{(d/r)^2}{s_1 s_2} 
\end{equation}
Also $w_1 + w_2 \le 2\frac{(d/r)}{\min(s_1,s_2)}$. Hence, for fixed
  $(s_1,s_2)$, 
\begin{equation}
\label{eq:easy:bound:fixed:s1s2}
\sum_{s_1 w_1 + s_2 w_2 = d/r} w_1 w_2 (w_1 + w_2) \min(s_1,s_2) \le 2
\frac{(d/r)^4}{s_1^2 s_2^2}
\end{equation}
(where we have used (\ref{eq:lambda:range}) to bound the number of
terms in the sum). Hence, the contribution to (\ref{eq:asymp:Nd:1:1})
from the terms with $s_1s_2\geq (d/r)/\log (d/r)$ is $o(d/r)^4$. 
(We  could use $(d/r)/f(d/r)$ for any function $f(d/r)$ which
is $o(d/r)$.)

\relax From now on we assume $(s_1,s_2)$ satisfy 
\begin{displaymath}
s_1s_2\leq (d/r)/\log (d/r).
\end{displaymath}

This implies that   
 $d/r-s_1w_1^0=o(d/r)$, and so  
the  sum (\ref{eq:asymp:Nd:1:1}) over such pairs $(s_1,s_2)$ becomes 
\begin{displaymath}
\sum_{r|d}\sum_{\stackrel{(s_1,s_2)=1}{s_1s_2\leq (d/r)/\log
(d/r)}}
\sum_{\lambda=0}^{\frac{d/r}{s_1s_2}}
r(\frac{d/r}{s_1}-\lambda s_2)(\lambda
s_1)(\frac{d/r}{s_1}-\lambda
s_2+\lambda s_1)\min (s_1,s_2)+o((d/r)^4)
\end{displaymath}

By replacing the inner sum with an integral over $\lambda$, the
sum becomes
\begin{equation}
\label{eq:covers} 
 \sum_{r|d}\frac{(d/r)^4}{12}
 \sum_{\stackrel{(s_1,s_2)=1}{s_1s_2\leq (d/r)/\log
(d/r)}}r(\frac{1}{s_1^3s_2^2}+
\frac{1}{s_1^2s_2^3})\min (s_1,s_2)+o((d/r)^4)
\end{equation}

We see immediately that for $r=1$, the sum is of the order $d^4$,
justifying the assumption that $r\leq d^{4/5}$.
  
Now we wish to further evaluate the sum (\ref{eq:covers}).  
 Recall,  
by definition, that the multiple $\zeta$ function
$\zeta(m_1,m_2)$
is defined by 
\begin{displaymath}
\zeta(m_1,m_2)=\sum_{s_1<s_2}^\infty
\frac{1}{s_1^{m_1}s_2^{m_2}}
\end{displaymath}
 
We will need to compute  $\zeta(2,2)$ and $\zeta(1,3)$. Such formulas
are well known (see e.g. \cite{Goncharov:multiple:zeta}), but we give
a short proof here for completeness. 
First we have
\begin{equation}
\label{eq:zeta} 
\zeta(2)^2=2\zeta(2,2)+\zeta(4)
\end{equation}
We also have 
\begin{multline*}
\zeta(2)^2=\sum_{x,y}
\frac{(x+y)^3}{x^2y^2(x+y)^3}=\sum\frac{(x^2+y^2)(x+y)+2x^2y+
2xy^2}{x^2y^2(x+y)^3}
= \\
\sum\left(\frac{1}{x^2}+\frac{1}{y^2}\right)\frac{1}{(x+y)^2}+2
\sum\left(\frac{1}{y}+\frac{1}{x}
\right)\frac{1}{(x+y)^3}=2\zeta(2,2)+4\zeta(1,3).
\end{multline*}
Together with (\ref{eq:zeta}) this gives 
\begin{displaymath}
\zeta(2,2)+\zeta(1,3)=\frac{\zeta(2)^2}{2}-
\frac{\zeta(4)}{4}=\zeta(4)
\end{displaymath}

 Thus 
 \begin{displaymath}
\sum_{s_1,s_2}(\frac{1}{s_1^3s_2^2}+\frac{1}{s_1^2s_2^3})
\min(s_1,s_2)=2(\zeta(2,2)+\zeta(1,3)+\zeta(4))=4\zeta(4)
\end{displaymath}
If we therefore perform the same sum conditioned over $(s_1,s_2)$
we introduce a factor of $\frac{1}{\zeta(4)}$ and so the sum is
$4$. We therefore find from (\ref{eq:covers}) that as $d\to\infty$, 
  
\begin{displaymath}N_d(1,1)=\sum_{r|d}
\frac{d^4}{3r^3}+o(d^4)
\end{displaymath}
 and so 
\begin{displaymath}
 \nu(\cH(1,1))=\lim_{D\to\infty}
\frac{10}{D^5}\sum_{d=1}^D (\sum_{r|d}
\frac{d^4}{3r^3}
+o(d^4))
\end{displaymath}

If we let $d=rq$, the  double sum is  evaluated by 
\begin{displaymath}
\sum_{r\leq D}\sum_{q\leq D/r} rq^4+o(D^5)=\sum_{r\leq
D}\frac{r}{5}\left(\frac{D}{r}\right)^5
+ro\left(\frac{D}{r}\right)^5+o(D^5)=\frac{\zeta(4)}{5}D^5+o(D^5)
\end{displaymath}

and so 
\begin{displaymath}
\nu(\cH(1,1))=\frac{\pi^4}{135}
\end{displaymath}

\bold{Proof of Lemma~\ref{lemma:volume:H2}.}
The number of covers of degree $d$ with $2$ cylinders is 
\begin{displaymath}
\sum_{\stackrel{h_1w_1+h_2w_2=d}{w_1 < w_2}} w_1w_2=
\sum_{r|d}
\sum_{\stackrel{(s_1,s_2)=1}{s_1w_1+s_2w_2=\frac{d}{r}}}
\frac{1}{2}w_1w_2
\end{displaymath}
Arguing as before, we can assume $d/r\to\infty$ as $d\to\infty$,
$s_1s_2\leq \frac{d/r} {\log(d/r)}$ and  find solutions of the
form 
$w_1=w_1^0-\lambda s_2$, $w_2=w_2^0+\lambda s_1$,
where 
 $\frac{s_1w_1^0}{(d/r)}\to 1$ as $d\to\infty$ and
$0\leq\lambda\leq \frac{d}{rs_1s_2}$.  
Thus the above sum is 
\begin{displaymath}
\frac{1}{2}\sum_{r|d}(\sum_{(s_1,s_2)=1}
\sum_{\lambda=0}^{\frac{d}{rs_1s_2}}
(\frac{d}{rs_1}-\lambda s_2)\lambda s_1+(o(d/r))^3)
\end{displaymath}
which if we replace by an integral over $\lambda$  becomes 
\begin{displaymath}
\frac{1}{2}\sum_{r|d}\sum_{(s_1,s_2)=1}
\frac{1}{s_1^2s_2^2}\frac{d^3}{r^3}
+o(d^3)
\end{displaymath}
which as we let $d\to\infty$ becomes
\begin{displaymath}
\frac{1}{12}
\sum_{r|d}\frac{d^3}{r^3}\frac{\zeta(2)^2}{\zeta(4)}+o(d^3)
\end{displaymath}

Up to lower order terms,  
the number of covers of degree $d$ with one
cylinder is  
\begin{displaymath}
\frac{1}{3}
\sum_{h|d}\sum_{l_1+l_2+l_3=
\frac{d}{h}}\frac{d}{h}=\sum_{r|d}\frac{d^3}{6r^3}
+o(d^3)
\end{displaymath}
Thus 
\begin{displaymath}
\nu(\cH(2))=\lim_{D\to\infty}\frac{8}{D^4}\sum_{d=1}^D\sum_{r|d}
\frac{d^3}{r^3}(\frac{\zeta(2)^2}{12\zeta(4)}+\frac{1}{6})
\end{displaymath}
which is equal to 
\begin{displaymath}
2\zeta(4)(\frac{\zeta(2)^2}{12\zeta(4)}+\frac{1}{6})
=\frac{\zeta(2)^2}{6}+\frac{\zeta(4)}{3}=\frac{\pi^4}{120}
\end{displaymath}

%--------------------------------

\bold{Proof of Lemma~\ref{lemma:asymp:NdP11}.}
We compute the asymptotics as $d \to \infty$ of
(\ref{eq:H11:fixed:r}):
Exactly as in Lemma~\ref{lemma:volume:H11} the first term is
$\frac{d^4}{3r^2}+o(d)^4$.  The second term is obviously
$o(d^4)$.

\bold{Proof of Lemma~\ref{lemma:asymp:NdP:2}.}
The sums here are exactly those found in
Lemma~\ref{lemma:volume:H2}. Thus the first in
(\ref{eq:NP:2}),
corresponding to the $2$-cylinder covers is asymptotic to
\begin{displaymath}
\frac{d^3}{12}
  \frac{\zeta(2)^2}
{\zeta(4)}\sum_{r|d}\frac{\mu(r)}{r^2} + o(d^3) = \frac{5}{24}
d^3
\sum_{r|d} \frac{\mu(r)}{r^2} + o(d^3)
\end{displaymath}
The last term in (\ref{eq:NP:2}), corresponding to the
$1$-cylinder
covers with a symmetry, is clearly $o(d^3)$. Finally, the second
term
in (\ref{eq:NP:2}), corresponding to $1$-cylinder covers without
a
symmetry is asymptotic to 
\begin{displaymath}
\frac{1}{6}\sum_{r|d}\mu(r)d(d/r)^2+ o(d^3) =
\frac{1}{6}d^3\sum_{r|d}\frac{\mu(r)}{r^2}+ o(d^3).
\end{displaymath}
This implies the lemma.
\qed\medskip

\bold{Proof of Theorem~\ref{theorem:constants:converge} for
$c(d)$.} 
By taking $s_1,u_2$ (or $s_2,u_1$) small, it is clear that the
inner sum in 
(\ref{eq:numerator}), is of order  $(d/r)^3$ as $d\to\infty$. 
 However the sum over both $u_1,u_2$ small has strictly smaller
order, so 
the contribution from the third term is of lower order than the
contribution from the first two; that is, of order $o((d/r)^3)$. 
Thus  by symmetry,
\begin{displaymath}
\int_{P_d(1,1)}\hat fd\nu=2\pi\epsilon^2\sum_{r|d}
\mu(r)\sum_{\stackrel{(s_1,s_2)=1}
{s_1u_1+s_2u_2=d/r}}\frac{u_1u_2(u_1+u_2)}{u_2^2}
\min(s_1,s_2)+o((d/r)^3)
\end{displaymath}
The inner sum is dominated by terms for which $u_1,s_2$ are large and $s_1,u_2$ are 
small. In particular, we can assume $s_1=\min(s_1,s_2)$ and $u_1+u_2=u_1$ up to terms of lower order.   We perform Mobius inversion on pairs $(s_1,s_2)$ so that the integral 
becomes 
\begin{displaymath}
2\pi\epsilon^2\sum_{r|d}\sum_{k|\frac{d}{r}}
\mu(r)\mu(k)
\sum_{s_1u_1+s_2u_2=\frac{d}{rk}}\frac{u_1^2}{u_2}ks_1
+o((d/r)^3)=
\end{displaymath}   
\begin{displaymath}
2\pi\epsilon^2\sum_{r|d}\sum_{k|\frac{d}{r}}\sum_{t|\frac{d}{rk}}
\mu(r)\mu(k)\sum_{\stackrel{(s_1,u_2)=1}
{s_1u_1+s_2u_2=\frac{d}{rkt}}}\frac{u_1^2}{u_2}ks_1
+o((d/r)^3)
\end{displaymath}
Replacing the sum with an integral as above, we find that
\begin{displaymath}
\int_{P_d(1,1)}\hat fd\nu
=2\pi\epsilon^2\sum_{r|d}\sum_{k|\frac{d}{r}}\sum_{t|\frac{d}{rk}}
\mu(r)\mu(k)\sum_{(s_1,u_2)=1}
\frac{d^3}
{3k^2r^3t^3s_1^2u_2^2}+
o((d/r)^3) 
\end{displaymath}
Now we set $y=rt$
and find
\begin{displaymath} 
\int_{P_d(1,1)}\hat fd\nu=\frac{2d^3\pi\epsilon^2}{3}
\sum_{(s_1,u_2)=1}\sum_{r|y|\frac{d}{k}}
\mu(r)\mu(k)
\frac{1}{k^2y^3s_1^2u_2^2}+o((d/r)^3)
\end{displaymath}
Now 
$\sum_{r|y}\mu(r)=0$ unless $y=1$, in which case it is $1$.
Thus 
\begin{displaymath}
\int_{P_d(1,1)}\hat fd\nu=\frac{2d^3
\pi\epsilon^2}{3}\sum_{k|d}\sum_{(s_1,u_2)=1}\frac{\mu(k)}{k^2}
\frac{1}{s_1^2u_2^2}
+o(d^3)
\end{displaymath}
Thus by (\ref{eq:c:d:immed}) and Lemma~\ref{lemma:asymp:NdP11} we find that 
\begin{displaymath}
\lim_{d\to\infty} c(d)=\lim_{d\to\infty}
\frac{\int_{P_d(1,1)}\hat fd\nu}{\nu(P_d(1,1))}
=5
\end{displaymath}
\qed\medskip

\bold{Proof of Theorem~\ref{theorem:constants:converge} for
$s_2(d)$.}
It is easy to see that only the first subcase of case 1 contributes as
$d \to \infty$. 

In view of Proposition~\ref{prop:value:s2:d},
the computation of the asymptotics 
in Lemma~\ref{lemma:asymp:NdP:2}, 
for the number of 
$2$  cylinder covers in $\cH(2)$ and 
Lemma~\ref{lemma:asymp:NdP11}, and   
\S\ref{sec:generic:constants}, 
\begin{displaymath}
\lim_{d\to\infty}s_2(d)=
 (\frac{5/24}{1/3})=\frac{5}{8}=s_2(1,1).
\end{displaymath}
\qed


\begin{thebibliography}{EMS2}

\bibitem{BO}
S.~Bloch and A.~Okounkov,
\newblock{``The Character of the Infinite Wedge Representation''},
\newblock{\em Adv. Math.}{\bf 149} (2000), no. 1, 1--60.

\bibitem{Dijkgraaf:texel}
R.~Dijkgraaf,
\emph{Mirror symmetry and elliptic curves},
The Moduli Space of Curves, R.~Dijkgraaf,
C.~Faber, G.~van~der~Geer (editors), 
Progress in Mathematics, \textbf{129},
Birkh\"auser, 1995.


\bibitem{Eskin:Masur:ae}
A.~Eskin, H.~Masur.
\newblock{``Pointwise asymptotic formulas on flat surfaces.''}
\newblock{\em Ergodic Theory Dynam. Systems} {\bf 21} (2001), no. 2, 443--478.

\bibitem{Eskin:Okounkov:volumes}
A.~Eskin, A.~Okounkov.
\newblock{``Asymptotics of numbers of branched coverings of a
torus
and volumes  of moduli spaces of holomorphic differentials''}
\newblock{\em To appear in Invent. Math.}

\bibitem{EOZ}
A.~Eskin, A.~Okounkov, A.~Zorich,
\newblock{``Volumes of spaces of Abelian and quadratic
  differentials''}
\newblock{\em In preparation.}


\bibitem{Eskin:Masur:Zorich:SV}
A.~Eskin, H.~Masur, A.~Zorich,
\emph{The Siegel-Veech constants},
in preparation.


\bibitem{Dani:Margulis:distribution}
S.G. Dani and G.A. Margulis.
\newblock {Limit distributions of orbits of unipotent flows and
values of
  quadratic forms}.
\newblock {\em Advances in Soviet Math.} {\bf 16}(1993), 91--137.


\bibitem{Goncharov:multiple:zeta}
A.~B.~Goncharov.
\newblock{Multiple polylogarithms, cyclotomy and modular complexes}. 
\newblock {\em Math. Res. Lett.} {\bf 5} (1998), no. 4, 497--516.

\bibitem{Gutkin:Judge:private}
E.~Gutkin, C.~Judge. 
\newblock{Private communication.}

\bibitem {G-J2} E. Gutkin, C. Judge. Affine mappings of
translation surfaces: geometry and arithmetic. Duke Math. J. {\bf 103} 
(2000) pp. 191-213.


\bibitem {Jacobson:rings} 
N.~Jacobson,
\newblock{\em The theory of rings.}
\newblock{Amer. Math. Soc. Colloq. Publ.} {\bf 37}, Amer. Math. Soc., 
Providence, RI, 1943. 


\bibitem{Shah:SL2}
N.~Shah.
\newblock {Limit distributions of expanding translates of certain
orbits on homogeneous spaces.}
\newblock {\em Proc. Indian Acad. Sci. (Math Sci)} {\bf 106}(2),
(1996), pp. 105--125.



\bibitem {Vee} W. Veech. Teichmuller curves in moduli space.
Eisenstein series and and an application to triangular billiards.
Invent. Math., 97, (1990), 117-171.



\bibitem{Veech:Siegel}
W. Veech, {\em Siegel measures},  {\em Annals of
Mathematics} {\bf 148} (1998), 895-944

\bibitem{Zemlyakov:Katok}
A.N.~Zemlyakov, A.B.~Katok.
\newblock{Topological transitivity of billiards in  
polygons.}
\newblock{\em Matem. Zametki} {\bf 18}(2) (1975) pp.~291--300.
English translation in {\em Math. Notes} {\bf 18}(2) (1976)  
pp.~760--764. 

\bibitem{Zorich:newton}
A.~Zorich. 
\newblock{``Square-Tiled Surfaces and Teichmuller Volumes of the
  Moduli Spaces of Abelian Differentials.''} 
\newblock{ \em To appear in the Newton Institute proceedings volume.}
\end{thebibliography}
\end{document}